\documentclass[oneside, 12pt]{amsart}
\usepackage{amsfonts}
\usepackage{amssymb}
\usepackage{latexsym}
\usepackage{graphics}
\usepackage[2emode]{psfrag}
\usepackage{mathrsfs}
\usepackage{amsthm}
\usepackage{amsmath}
\usepackage[all]{xy}
\usepackage{enumerate}
\usepackage{hyperref}

\begin{document}

\newcommand{\thlabel}[1]{\label{th:#1}}
\newcommand{\thref}[1]{Theorem~\ref{th:#1}}
\newcommand{\selabel}[1]{\label{se:#1}}
\newcommand{\seref}[1]{Section~\ref{se:#1}}
\newcommand{\lelabel}[1]{\label{le:#1}}
\newcommand{\leref}[1]{Lemma~\ref{le:#1}}
\newcommand{\prlabel}[1]{\label{pr:#1}}
\newcommand{\prref}[1]{Proposition~\ref{pr:#1}}
\newcommand{\colabel}[1]{\label{co:#1}}
\newcommand{\coref}[1]{Corollary~\ref{co:#1}}
\newcommand{\relabel}[1]{\label{re:#1}}
\newcommand{\reref}[1]{Remark~\ref{re:#1}}
\newcommand{\exlabel}[1]{\label{ex:#1}}
\newcommand{\exref}[1]{Example~\ref{ex:#1}}
\newcommand{\delabel}[1]{\label{de:#1}}
\newcommand{\deref}[1]{Definition~\ref{de:#1}}
\newcommand{\eqlabel}[1]{\label{eq:#1}}
\newcommand{\equref}[1]{(\ref{eq:#1})}
\newcommand{\norm}[1]{\| #1 \|}

\def\N{\mathbb N}
\def\Z{\mathbb Z}
\def\Q{\mathbb Q}
\def\mod{\textit{\emph{~mod~}}}
\def\R{\mathcal R}
\def\S{\mathcal S}
\def\*C{{^*\mathcal C}}
\def\C{\mathcal C}
\def\D{\mathcal D}
\def\J{\mathcal J}
\def\M{\mathcal M}
\def\T{\mathcal T}

\newcommand{\Hom}{{\rm Hom}}
\newcommand{\End}{{\rm End}}
\newcommand{\Fun}{{\rm Fun}}
\newcommand{\Mor}{{\rm Mor}\,}
\newcommand{\Aut}{{\rm Aut}\,}
\newcommand{\Hopf}{{\rm Hopf}\,}
\newcommand{\Ann}{{\rm Ann}\,}
\newcommand{\Ker}{{\rm Ker}\,}
\newcommand{\Coker}{{\rm Coker}\,}
\newcommand{\im}{{\rm Im}\,}
\newcommand{\coim}{{\rm Coim}\,}
\newcommand{\Trace}{{\rm Trace}\,}
\newcommand{\Char}{{\rm Char}\,}
\newcommand{\Mod}{{\bf mod}}
\newcommand{\Spec}{{\rm Spec}\,}
\newcommand{\Span}{{\rm Span}\,}
\newcommand{\sgn}{{\rm sgn}\,}
\newcommand{\Id}{{\rm Id}\,}
\newcommand{\Com}{{\rm Com}\,}
\newcommand{\codim}{{\rm codim}}
\newcommand{\Mat}{{\rm Mat}}
\newcommand{\can}{{\rm can}}
\newcommand{\sign}{{\rm sign}}
\newcommand{\kar}{{\rm kar}}
\newcommand{\rad}{{\rm rad}}

\def\Ab{\underline{\underline{\rm Ab}}}
\def\lan{\langle}
\def\ran{\rangle}
\def\ot{\otimes}

\def\id{\textrm{{\small 1}\normalsize\!\!1}}
\def\To{{\multimap\!\to}}
\def\bigperp{{\LARGE\textrm{$\perp$}}}
\newcommand{\QED}{\hspace{\stretch{1}}
\makebox[0mm][r]{$\Box$}\\}

\def\AA{{\mathbb A}}
\def\BB{{\mathbb B}}
\def\CC{{\mathbb C}}
\def\DD{{\mathbb D}}
\def\EE{{\mathbb E}}
\def\FF{{\mathbb F}}
\def\GG{{\mathbb G}}
\def\HH{{\mathbb H}}
\def\II{{\mathbb I}}
\def\JJ{{\mathbb J}}
\def\KK{{\mathbb K}}
\def\LL{{\mathbb L}}
\def\MM{{\mathbb M}}
\def\NN{{\mathbb N}}
\def\OO{{\mathbb O}}
\def\PP{{\mathbb P}}
\def\QQ{{\mathbb Q}}
\def\RR{{\mathbb R}}
\def\SS{{\mathbb S}}
\def\TT{{\mathbb T}}
\def\UU{{\mathbb U}}
\def\VV{{\mathbb V}}
\def\WW{{\mathbb W}}
\def\XX{{\mathbb X}}
\def\YY{{\mathbb Y}}
\def\ZZ{{\mathbb Z}}

\def\aa{{\mathfrak A}}
\def\bb{{\mathfrak B}}
\def\cc{{\mathfrak C}}
\def\dd{{\mathfrak D}}
\def\ee{{\mathfrak E}}
\def\ff{{\mathfrak F}}
\def\gg{{\mathfrak G}}
\def\hh{{\mathfrak H}}
\def\ii{{\mathfrak I}}
\def\jj{{\mathfrak J}}
\def\kk{{\mathfrak K}}
\def\ll{{\mathfrak L}}
\def\mm{{\mathfrak M}}
\def\nn{{\mathfrak N}}
\def\oo{{\mathfrak O}}
\def\pp{{\mathfrak P}}
\def\qq{{\mathfrak Q}}
\def\rr{{\mathfrak R}}
\def\ss{{\mathfrak S}}
\def\tt{{\mathfrak T}}
\def\uu{{\mathfrak U}}
\def\vv{{\mathfrak V}}
\def\ww{{\mathfrak W}}
\def\xx{{\mathfrak X}}
\def\yy{{\mathfrak Y}}
\def\zz{{\mathfrak Z}}

\def\aaa{{\mathfrak a}}
\def\bbb{{\mathfrak b}}
\def\ccc{{\mathfrak c}}
\def\ddd{{\mathfrak d}}
\def\eee{{\mathfrak e}}
\def\fff{{\mathfrak f}}
\def\ggg{{\mathfrak g}}
\def\hhh{{\mathfrak h}}
\def\iii{{\mathfrak i}}
\def\jjj{{\mathfrak j}}
\def\kkk{{\mathfrak k}}
\def\lll{{\mathfrak l}}
\def\mmm{{\mathfrak m}}
\def\nnn{{\mathfrak n}}
\def\ooo{{\mathfrak o}}
\def\ppp{{\mathfrak p}}
\def\qqq{{\mathfrak q}}
\def\rrr{{\mathfrak r}}
\def\sss{{\mathfrak s}}
\def\ttt{{\mathfrak t}}
\def\uuu{{\mathfrak u}}
\def\vvv{{\mathfrak v}}
\def\www{{\mathfrak w}}
\def\xxx{{\mathfrak x}}
\def\yyy{{\mathfrak y}}
\def\zzz{{\mathfrak z}}

\newcommand{\aA}{\mathscr{A}}
\newcommand{\bB}{\mathscr{B}}
\newcommand{\cC}{\mathscr{C}}
\newcommand{\dD}{\mathscr{D}}
\newcommand{\eE}{\mathscr{E}}
\newcommand{\fF}{\mathscr{F}}
\newcommand{\gG}{\mathscr{G}}
\newcommand{\hH}{\mathscr{H}}
\newcommand{\iI}{\mathscr{I}}
\newcommand{\jJ}{\mathscr{J}}
\newcommand{\kK}{\mathscr{K}}
\newcommand{\lL}{\mathscr{L}}
\newcommand{\mM}{\mathscr{M}}
\newcommand{\nN}{\mathscr{N}}
\newcommand{\oO}{\mathscr{O}}
\newcommand{\pP}{\mathscr{P}}
\newcommand{\qQ}{\mathscr{Q}}
\newcommand{\rR}{\mathscr{R}}
\newcommand{\sS}{\mathscr{S}}
\newcommand{\tT}{\mathscr{T}}
\newcommand{\uU}{\mathscr{U}}
\newcommand{\vV}{\mathscr{V}}
\newcommand{\wW}{\mathscr{W}}
\newcommand{\xX}{\mathscr{X}}
\newcommand{\yY}{\mathscr{Y}}
\newcommand{\zZ}{\mathscr{Z}}

\newcommand{\Aa}{\mathcal{A}}
\newcommand{\Bb}{\mathcal{B}}
\newcommand{\Cc}{\mathcal{C}}
\newcommand{\Dd}{\mathcal{D}}
\newcommand{\Ee}{\mathcal{E}}
\newcommand{\Ff}{\mathcal{F}}
\newcommand{\Gg}{\mathcal{G}}
\newcommand{\Hh}{\mathcal{H}}
\newcommand{\Ii}{\mathcal{I}}
\newcommand{\Jj}{\mathcal{J}}
\newcommand{\Kk}{\mathcal{K}}
\newcommand{\Ll}{\mathcal{L}}
\newcommand{\Mm}{\mathcal{M}}
\newcommand{\Nn}{\mathcal{N}}
\newcommand{\Oo}{\mathcal{O}}
\newcommand{\Pp}{\mathcal{P}}
\newcommand{\Qq}{\mathcal{Q}}
\newcommand{\Rr}{\mathcal{R}}
\newcommand{\Ss}{\mathcal{S}}
\newcommand{\Tt}{\mathcal{T}}
\newcommand{\Uu}{\mathcal{U}}
\newcommand{\Vv}{\mathcal{V}}
\newcommand{\Ww}{\mathcal{W}}
\newcommand{\Xx}{\mathcal{X}}
\newcommand{\Yy}{\mathcal{Y}}
\newcommand{\Zz}{\mathcal{Z}}


\numberwithin{equation}{section}

\renewcommand{\theequation}{\thesection.\arabic{equation}}


\newcommand{\bara}[1]{\overline{#1}}
\newcommand{\ContFunt}[2]{\bara{\mathrm{Funt}}(#1,\,#2)}
\newcommand{\lBicomod}[2]{{}_{#1}\mM^{#2}}
\newcommand{\rBicomod}[2]{{}^{#1}\mM_{#2}}
\newcommand{\Bicomod}[2]{{}^{#1}\mM^{#2}}
\newcommand{\lrBicomod}[2]{{}^{#1}\mM^{#2}}
\newcommand{\rcomod}[1]{\mM^{#1}}
\newcommand{\rmod}[1]{\mM_{#1}}
\newcommand{\Bimod}[2]{{}_{#1}\mM_{#2}}
\newcommand{\Sf}[1]{\mathsf{#1}}
\newcommand{\Sof}[1]{S^{{\bf #1}}}
\newcommand{\Tof}[1]{T^{{\bf #1}}}
\renewcommand{\hom}[3]{\mathrm{Hom}_{#1}\left(\underset{}{}#2,\,#3\right)}
\newcommand{\Coint}[2]{\mathrm{Coint}(#1,#2)}
\newcommand{\IntCoint}[2]{\mathrm{InCoint}(#1,#2)}
\newcommand{\Coder}[2]{\mathrm{Coder}(#1,#2)}
\newcommand{\IntCoder}[2]{\mathrm{InCoder}(#1,#2)}
\newcommand{\lr}[1]{\left(\underset{}{} #1 \right)}
\newcommand{\Ext}[3]{\mathrm{Ext}_{\eE}^{#1}\lr{#2,\,#3}}
\newcommand{\equalizerk}[2]{\mathfrak{eq}_{#1,\,#2}^k}
\newcommand{\equalizer}[2]{\mathfrak{eq}_{#1,\,#2}}
\newcommand{\coring}[1]{\mathfrak{#1}}
\newcommand{\tensor}[1]{\otimes_{#1}}


\def\units{{\mathbb G}_m}
\def\rightact{\hbox{$\leftharpoonup$}}
\def\leftact{\hbox{$\rightharpoonup$}}

\def\*C{{}^*\hspace*{-1pt}{\Cc}}

\def\text#1{{\rm {\rm #1}}}

\def\smashco{\mathrel>\joinrel\mathrel\triangleleft}

\def\cosmash{\mathrel\triangleright\joinrel\mathrel<}

\def\ol{\overline}
\def\ul{\underline}
\def\dul#1{\underline{\underline{#1}}}
\def\Nat{\dul{\rm Nat}}
\def\Set{\dul{\rm Set}}

\renewcommand{\subjclassname}{\textup{2000} Mathematics Subject Classification}

\newtheorem{prop}{Proposition}[section]
\newtheorem{lemma}[prop]{Lemma}
\newtheorem{cor}[prop]{Corollary}
\newtheorem{theo}[prop]{Theorem}

\theoremstyle{definition}
\newtheorem{Def}[prop]{Definition}
\newtheorem{ex}[prop]{Example}
\newtheorem{exs}[prop]{Examples}

\theoremstyle{remark}
\newtheorem{rems}[prop]{Remarks}
\newtheorem{rem}[prop]{Remark}

\title[Cohomology for Bicomodules. separable and Maschke functors]
{Cohomology for Bicomodules. Separable and Maschke functors.}
\date{August 2006}
\author{L. El Kaoutit}
\address{Departamento de \'Algebra. Facultad de Educaci\'on y Humanidades de Ceuta.
Universidad de Granada. El Greco N 10. E51002 Ceuta, Spain}
\email{kaoutit@ugr.es}


\author{J. Vercruysse}
\address{Faculty of Applied Sciences, Vrije Universiteit Brussel (VUB), B-1050 Brussels, Belgium}
\email{jvercruy@vub.ac.be}
\urladdr{homepages.vub.ac.be/\~{}jvercruy}


\baselineskip 16pt

\begin{abstract}
We introduce the category of bicomodules for a comonad in a
Grothendieck category whose underlying functor is right exact and
preserves direct sums. We characterize comonads with a separable
forgetful functor by means of cohomology groups using cointegrations
into bicomodules. We present two applications: the characterization
of coseparable corings stated in \cite{Guzman:1989}, and the
characterization of coseparable coalgebras coextensions stated in
\cite{N}.
\end{abstract}

\maketitle

\section*{Introduction}
In \cite{Jonah:1968} D. W. Jonah studied the second and the third
cohomology groups of coalgebras defined in a, not necessary abelian,
multiplicative category (see also
\cite{Ardizzoni/Menini/Dragos:2006}). M. Kleiner gave in
\cite{Kleiner:1985} a cohomological characterization of separable
algebras using integrations. Another approach via derivations was
given by M. Barr and G. Rinehart in \cite{BR}. This last one has
been dualised to the case of coseparable coalgebras by Doi
\cite{Doi}. Nakajima \cite{N} showed that Doi's results can be
extended to the coalgebra extensions (or co-extension) with a
co-commutative base coalgebra. In \cite{Guzman:1989}, F. Guzman used
Jonah's methods to generalize Doi's characterization for corings
over an arbitrary base-ring and unified this with a dualisation of
Kleiner's approach of cointegrations. This gives rise to a nice
characterization of coseparable corings in terms of cohomology,
derived functors and both cointegrations and coderivations.
Unfortunately  this last characterization can not be applied to
coalgebra co-extensions, and Nakajima's results is not recovered.

The common framework behind Guzman's and Nakajima's approach is the
fact that both coseparable corings and coseparable coalgebra
co-extensions can be interpret as comonads with a separable
forgetful functor (in the sense of
\cite{Nastasescu/Bergh/Oystaeyen:1989}, see below). In all
situations discussed before, the multiplicative base-category was
additive with cokernels and arbitrary direct sums, and the (co)monad
functor preserved cokernels and direct sums. In the present paper we
will approach the problem by this comonad point of view. We work
with a comonad over a Grothendieck category (not necessary
multiplicative) whose underlying functor fits the above mentioned
class of functors. These functors were studied in relation with
corings in \cite{Gomez-Torrecillas:2002}, see also
\cite{Kaoutit:2006} and references sated there. We will present a
generalisation of Guzman's characterization in this situation, and
as a particular application we also give, under different
assumption, Nakajima's result.

We will start by defining the category of bicomodule over this
comonad as in \cite{Jonah:1968}, and we consider its universal
cogenerator \cite{Eilenberg/Moore:1965} (i.e. the universal
adjunction defining the comonadic structure) in order to prove that
the forgetful functor in this universal adjunction is separable
(\cite{Nastasescu/Bergh/Oystaeyen:1989}, see below) if and only if
the forgetful functor in bicomodules is Maschke
(\cite{Caenepeel/Militaru:2003}, see below) if and only if the
comultiplication splits in the category of bicomodules. This will be
the main result of section \ref{Sect-1} (Theorem \ref{Separable})
(see \cite{Ard} for a different approach). In section \ref{Sect-2}
we define cointegrations and coderivations, we also establish, as in
\cite{Guzman:1989}, an isomorphism between the abelian group of
cointegrations into a comonad and the group of all coderivations.
This will serve to show that the comultiplication splits as a
morphism of bicomodules if and only if the universal cointegration
is inner if and only if the universal coderivation is inner
(Corollary \ref{split-sequence}). Section \ref{Sect-3} is devoted to
the relative cohomology for bicomodules defined as in
\cite{Jonah:1968} using a relative resolution with respect to the
injective class of sequences in the category of bicomodules which
are cosplit after forgetting the left coaction. Up to isomorphisms,
cointegrations appear as $1$-cocycles and inner cointegrations as
$1$-coboundaries. The relative injectivity is thus interpreted by
the fact that all into-cointegrations are inner. This happens for
all bicomodules if and only if the comultiplication splits in the
category of bicomodules (Theorem \ref{characterization}). The last
section presents two applications of this last theorem, the first
one makes use of the comonad defined by tensor product over algebras
\cite{Guzman:1989}, and the second uses cotensor product over
coalgebras over fields \cite{N}. Since in recent years it became
clear that corings and comodules provide a general framework to
study entwining structures and entwined modules and by this all
sorts of (relative) Hopf-modules (we refer to
\cite{Brzezinski/Wisbauer:2003} for a profound overview),
separability properties of these structures are covered by our
theory as special cases.

\bigskip
{\textsc{Notations and Basic Notions}:} Given any Hom-set
category $\aA$, the notation $X \in \aA$ means that $X$ is an object
of $\aA$. The identity morphism of $X$ will be denoted by $X$
itself. The set of all morphisms $f:X \to X'$ in $\aA$, is denoted
by $\hom{\aA}{X}{X'}$. The identity functor of $\aA$ will be denoted
by $\id_{\aA}: \aA \rightarrow \aA$. A natural transformation
between two functor $\Ff,\,\Gg: \aA \to \bB$, is denoted by
$\beta_{-}: \Ff \to \Gg$. If $\Hh: \bB \to \cC$, and $\Ii: \dD \to
\aA$ are other functors. Then, $\beta_{\Ii(-)}$ (or $\beta_\Ii$)
denotes the natural transformation defined at each object $Z \in
\dD$ by $\beta_{\Ii(Z)}:\Ff\Ii(Z) \to \Gg\Ii(Z)$, while $\Hh
\beta_{-}$ (or $\Hh\beta$) denotes the natural transformation
defined at each object $X \in \aA$ by $\Hh(\beta_X): \Hh\Ff(X) \to \Hh\Gg(X)$.\\
Any covariant functor $\Ff:\aA \to \bB$ leads to a (bi)functor
$$\Hom_\bB(\Ff(-),\Ff(-)):\aA^{op}\times\aA\to \sS et.$$
In particular, the identical functor $\id_\aA:\aA\to\aA$ gives rise
to $$\Hom_\aA(-,-):\aA^{op}\times \aA \to \sS et.$$ So we find a
natural transformation induced by $\Ff$,
$$\fF:\Hom_\aA(-,-)\to\Hom_\bB(\Ff(-),\Ff(-));$$
defined by $\fF_{X,\,X'}(f)=\Ff(f)$, for any arrow  $f:X\to X'$ in
$\aA$. Recall from \cite{Nastasescu/Bergh/Oystaeyen:1989} that the
functor $\Ff$ is called \emph{separable} if and only if $\fF$ has a
left inverse, i.e. there exists a natural transformation
$$\pP:\Hom_\bB(\Ff(-), \Ff(-))\to\Hom_\aA(-,-)$$
such that $\pP\circ\fF=\id_{\Hom_\aA(-,-)}$. If in addition $\Ff$
has a right adjoint functor $\Gg: \bB \to \aA$ with unit $\eta_{-}:
\id_\aA \to \Gg \Ff$. Then, it is well known from
\cite{Rafael:1990}, that $\Ff$ is separable if and only if there
exists a natural transformation $\mu: \Gg \Ff \to \id_\aA$ such that
$\mu \circ \eta \,=\, \id_\aA$.\\ Let $\Ff: \aA\to\bB$ be again a
covariant functor. Recall from \cite{Caenepeel/Militaru:2003}, that
an object $M\in\aA$ is called \emph{relative injective} (or
\emph{$\Ff$-injective}) if and only if for every morphism $i:X\to
X'$ in $\aA$, such that $\Ff(i):\Ff(X)\to \Ff(X')$ has a left
inverse $j$ in $\Bb$ (i.e. $\Ff(i)$ is a split monomorphism or just
split-mono) and for every $f:X\to M$ in $\aA$ we can find a morphism
$g: X'\to M$ in $\aA$ such that $g\circ i=f$. The functor $\Ff$ is
said to be \emph{a Maschke functor} if every object of $\aA$ is
$\Ff$-injective. If in addition  $\Ff$ has a right adjoint functor
$\Gg: \bB \to \aA$ with unit $\eta_{-}:\id_\aA \to \Gg \Ff$. Then,
by \cite[Theorem 3.4]{Caenepeel/Militaru:2003}, an object $M \in
\aA$ is $\Ff$-injective if and only if $\eta_{M}$  has a left
inverse. In particular $\Ff$ is a Maschke functor if and only if for
every object $M \in \aA$, $\eta_{M}$ has a left inverse.

Assume that a preadditive category $\aA$ is given. Following to
\cite{Jonah:1968}, a sequence $$\xymatrix{E: \, X \ar@{->}^-{i}[r] &
X \ar@{->}^-{j}[r] & X''}$$ (i.e. $j \circ i=0$) is said to be
\emph{co-exact} if $i$ has a cokernel and if in the commutative
diagram
$$\xymatrix@C=40pt{
X \ar@{->}^-{i}[r] & X' \ar@{->}^-{j}[r] \ar@{->}_-{i^c}[d] & X'' \\
& {\rm Coker}(i) \ar@{-->}_-{l}[ur] & }$$ $l$ is a monomorphism. If
in addition $l$ is a split-mono, then $E$ is said to be
\emph{cosplit}. The \emph{exact} and \emph{split} sequence are
dually defined by using kernels. The notations of sequences,
coexact, cosplit,... are extended to long diagrams simply by
applying them to each consecutive pair of morphisms. One can prove
that the above notions of exact and coexact sequences coincide with
the usual meaning of exact sequences in abelian categories. In case
of diagrams of the form
$$\xymatrix{E': \, 0 \ar@{->}[r] & X \ar@{->}[r] & X' \ar@{->}[r] &
X'' \ar@{->}[r] & 0}$$ (i.e. short sequence) in the category $\aA$ ,
we have by \cite[Lemma 2.1]{Jonah:1968}, that $E'$ is cosplit if and
only if it is split.

Let $\eE$ be a class of sequences in $\aA$, then an object $X \in
\aA$ is said to be \emph{$\eE$-injective} if $\Hom_{\aA}(E,X)$ is an
exact sequence of abelian groups, for every sequence $E$ in $\eE$.
The class of all $\eE$-injective objects is denoted by $\iI_\eE$.
Conversely, given $\iI$ a class of objects of $\aA$, a sequence $E$
of morphism of $\aA$ is said to be \emph{$\iI$-exact} if
$\Hom_{\aA}(E,Y)$ is an exact sequence of an abelian groups, for
every object $Y$ in $\iI$. The class of all $\iI$-exact sequences is
denoted by $\eE_\iI$. A class of sequences $\eE$ in $\aA$ is said to
be \emph{closed} whenever $\eE$ coincides with $\eE_{\iI_\eE}$. An
\emph{injective class} is a closed class of sequences $\eE$ such
that, for every morphism $X \to X'$, there exists a morphism $X' \to
Y$ with $Y \in \iI_\eE$ and with $X \to X' \to Y$ in $\eE$. If in
addition the category $\aA$ poses cokernels, then one can check that
the class $\eE_0$ of all cosplit sequences form an injective class
and $\iI_{\eE_0}$ is exactly the class of all objects of $\aA$.
Given any adjunction $\xymatrix{ \fF: \aA \ar@<0,5ex>[r] &
\ar@<0,5ex>[l] \bB: \gG}$ with $\fF$ is left adjoint functor to
$\gG$ (we use the notation $\fF \dashv \gG$), and a class of
sequences $\eE'$ in $\bB$. Denote by $\eE=\fF^{-1}(\eE')$ the class
of sequences $E$ in $\aA$ such that $\fF(E)$ is in $\eE'$. The
Eilenberg-Moore Theorem \cite[Theorem 2.9]{Jonah:1968} asserts that
$\eE$ is an injective class whenever $\eE'$ is.

\section{Bicomodules and Separability}\label{Sect-1}

Let $\Aa$ and $\Bb$ two Grothendieck categories, we denote by
$\ContFunt{\Aa}{\Bb}$ the class of all (additive) covariant functors
$F: \Aa \rightarrow \Bb$ such that $F$ preserves cokernels and
commutes with direct sums. Thus $F$ commutes with inductive limits.
By \cite[Lemma 5.1]{Iglesias/Gomez/Nastasescu:1999}, the natural
transformations between two objects of the class
$\ContFunt{\Aa}{\Bb}$ form a set. Henceforth, $\ContFunt{\Aa}{\Bb}$
is a Hom-set category (or Set-category).

A \emph{comonad} in a category $\Aa$ is a three-tuple ${\bf F}
=(F,\delta,\xi)$ consisting of an endo-functor $F:\Aa \rightarrow
\Aa$ and two natural transformations $\delta: F \rightarrow F^2=F
\circ F$ and $\xi:F \rightarrow \id_{\Aa}$  such that
\begin{equation}\label{Comonad}
\delta_F \, \circ \, \delta \,\, =\,\, F\delta \,\circ \, \delta,
\quad F\xi \,\circ \, \delta \,\, =\,\, \xi_F \,\circ \, \delta \,\,
= \,\, F,
\end{equation}
where we denote the identical natural transformation $F \rightarrow
F$ again by $F$.

It is well known from
\cite{Huber:1961,Eilenberg/Moore:1965,Kleisli:1965}, that any
adjunction $\xymatrix{ S: \Bb \ar@<0,5ex>[r] & \ar@<0,5ex>[l] \Aa:
T}$ with $S \dashv T$, leads to a comonad in $\Aa$ given by the
three-tuple $(ST, S\eta_{T}, \zeta)$, where $\eta: \id_{\Bb}
\rightarrow TS$ and $\zeta:ST \rightarrow \id_{\Aa}$ are,
respectively, the unit and the counit of this adjunction.

Let ${\bf F}=(F,\delta,\xi)$ be a comonad in $\Aa$ with $F \in
\ContFunt{\Aa}{\Aa}$. We define the category of $(\Bb,{\bf
F})$-\emph{bicomodules} $\lBicomod{\Bb}{{\bf F}}$ by the following
data:
\begin{enumerate}[$\bullet$]
\item \emph{Objects}: A $(\Bb,\bf F)$-bicomodule is a pair
$(\mathsf{M},\mmm)$ consisting of a functor $\Sf{M} \in
\ContFunt{\Bb}{\Aa}$ and natural transformation $\mmm: \Sf{M}
\rightarrow F\Sf{M}$ satisfying
\begin{equation}\label{comod-1}
\delta_{\Sf{M}} \,\circ \, \mmm \,\, =\,\, F\mmm \,\circ \,
\mmm,\quad \xi_{\Sf{M}} \, \circ \, \mmm \,\, = \,\, \Sf{M}.
\end{equation}

\item \emph{Morphisms}: A morphism $\Sf{f}: (\Sf{M},\mmm) \rightarrow
(\Sf{M}',\mmm')$ is a natural transformation $\Sf{f}:\Sf{M}
\rightarrow \Sf{M}'$ satisfying
\begin{equation}\label{comod-2}
\mmm' \,\circ \, \Sf{f} \,\, =\,\, F\Sf{f} \, \circ \, \mmm.
\end{equation}
\end{enumerate}

It is easily seen that $(F\Sf{M},\delta_{\Sf{M}})$ is an object of
the category $\lBicomod{\Bb}{{\bf F}}$, for every object $\Sf{M}
\in \ContFunt{\Bb}{\Aa}$. This in fact establishes a functors
$\fF: \ContFunt{\Bb}{\Aa} \rightarrow \lBicomod{\Bb}{{\bf F}}$
with a left adjoint the forgetful functor $\oO:
\lBicomod{\Bb}{{\bf F}} \rightarrow \ContFunt{\Bb}{\Aa}$.

Similarly, we can define the category of $({\bf
F},\Bb)$-bicomodules denoted by $\rBicomod{{\bf F}}{\Bb}$, using
this time the objects of the category $\ContFunt{\Aa}{\Bb}$.

\begin{rem}
Given any adjunction $\xymatrix{ \Sf{M}: \Bb \ar@<0,5ex>[r] &
\ar@<0,5ex>[l] \Aa: \Sf{N}}$ such that $\Sf{M} \dashv \Sf{N}$ with
counit $\zeta$ and unit $\eta$. Then \cite[Proposition
1.1]{Gomez-Torrecillas:2006} establishes an one-to-one
correspondences between natural transformations $\mmm: \Sf{M} \to
F\Sf{M}$ satisfying equation \eqref{comod-1} and homomorphisms of
comonads from $(\Sf{M}\Sf{N},\Sf{M} \eta_{\Sf{N}}, \zeta)$ to
$\bf{F}$, and a natural transformations $\sss: \Sf{N} \to \Sf{N}F$
satisfying the dual version of equation \eqref{comod-1}. When
$\Sf{N}$ and $\Sf{M}$ are both right exact and preserve direct sums,
then the previous correspondence can be interpreted in our
terminology as follows: There are bijections between the bicomodule
structures on $\Sf{M}$, the bicomodule structures on $\Sf{N}$, and
the homomorphisms of comonads from $(\Sf{M}\Sf{N},\Sf{M}
\eta_{\Sf{N}}, \zeta)$ to $\bf{F}$.
\end{rem}

Take now ${\bf G}=(G,\vartheta,\varsigma)$ another comonad in $\Bb$
with $G \in \ContFunt{\Bb}{\Bb}$, we define the category of $({\bf
G},{\bf F})$-\emph{bicomodules} $\lrBicomod{{\bf G}}{{\bf F}}$ as
follows:

\begin{enumerate}[$\bullet$]
\item \emph{Objects}: A $(\bf G,\bf F)$ bicomodule is a three-tuple $(\mathsf{M},\mmm,\nnn)$
consisting of a functor $\Sf{M} \in \ContFunt{\Bb}{\Aa}$ and two
natural transformations $\mmm: \Sf{M} \rightarrow F\Sf{M}$ , $\nnn:
\Sf{M} \rightarrow \Sf{M} G$ such that $(\Sf{M},\mmm) \in
\lBicomod{\Bb}{{\bf F}}$ and $(\Sf{M},\nnn) \in \rBicomod{{\bf
G}}{\Aa}$, that is
\begin{equation}\label{comod-00}
\delta_{\Sf{M}} \circ \mmm\,\, =\,\, F\mmm \circ \mmm,\,\,
\xi_{\Sf{M}} \circ \mmm \,\, =\,\, \Sf{M}\,\,\, \text{ and }\,\,\,
\Sf{M}\vartheta \circ \nnn \,\, =\,\, \nnn_G \circ \nnn, \,\,
\Sf{M}\varsigma \circ \nnn\,\,=\,\, \Sf{M}
\end{equation} with compatibility condition
\begin{equation}\label{comod-11}
\mmm_G \,\circ \, \nnn \,\, =\,\, F\nnn \,\circ \, \mmm.
\end{equation}In other words $\mmm$ is a morphism of $\rBicomod{{\bf
G}}{\Aa}$, equivalently, $\nnn$ is a morphism of
$\lBicomod{\Bb}{{\bf F}}$, where $(F\Sf{M},F\nnn) \in \rBicomod{{\bf
G}}{\Aa}$ and $(\Sf{M}G,\mmm_G) \in \lBicomod{\Bb}{{\bf F}}$.

\item \emph{Morphisms}: A morphism $\Sf{f}: (\Sf{M},\mmm,\nnn) \rightarrow
(\Sf{M}',\mmm',\nnn')$ is a natural transformation $\Sf{f}:\Sf{M}
\rightarrow \Sf{M}'$ such that $\Sf{f}: (\Sf{M},\mmm) \rightarrow
(\Sf{M}',\mmm')$ is a morphism of $\lBicomod{\Bb}{{\bf F}}$ and
$\Sf{f}: (\Sf{M},\nnn) \rightarrow (\Sf{M}',\nnn')$ is a morphism of
$\rBicomod{{\bf G}}{\Aa}$, that is
\begin{equation}\label{bicolineal}
\nnn' \circ \Sf{f}\,\, =\,\, \Sf{f}_G \circ \nnn\,\,\, \text{and}
\,\,\, \mmm' \circ \Sf{f} \,\, =\,\, F\Sf{f} \circ \mmm.
\end{equation}
\end{enumerate}

It is clear that $\lrBicomod{{\bf \id_{\Bb}}}{{\bf F}} =
\lBicomod{\Bb}{{\bf F}}$ and $\lrBicomod{{\bf
G}}{\id_{\Aa}}=\rBicomod{{\bf G}}{\Aa}$, where ${\bf \id_{\Aa}}$ and
${\bf \id_{\Bb}}$ are endowed with a trivial comonad structure.

\begin{rem}
It is easily seen that $\ContFunt{\Aa}{\Aa}$ is a strict monoidal
category (or multiplicative category), taking the composition of
functors as the tensor product  and $\id_\Aa$ as the unit object.
To any coalgebra in a monoidal category one can associate in a
canonical way a category of bicomodules, see \cite[Section
1]{Jonah:1968}. If we consider $\bf F$ as a coalgebra in
$\ContFunt{\Aa}{\Aa}$, then the category of $(\bf F,\bf
F)$-bicomodules as defined above coincides exactly with this
canonical one. However, if we consider $(\bf G, \bf
F)$-bicomodules and thus the base-category is changed, the
monoidal arguments fail. In that case one must consider the
$2$-category of Grothendieck categories ($0$-cells), additive
functors that preserve cokernels and commute with direct sums
($1$-cells), and natural transformations ($2$-cells). Observe that
$\bf F$ and $\bf G$ are comonads inside this $2$-category (see
\cite{Kaoutit:2006} for elementary treatment).
\end{rem}

By the observation that the bicomodules as introduced above
coincide with certains $1$-cells in a $2$-category, we can state
the following well known lemma.

\begin{lemma}\label{conucleos}
Let $\Aa$ (respectively $\Bb$) be a Grothendieck category, and
${\bf{F}} \,=\, (F,\delta,\xi)$ (respectively ${\bf{G}}\, =\,
(G,\vartheta,\varsigma)$) a comonad in $\Aa$ (respectively in $\Bb$)
whose underlying functor $F$ (respectively $G$) is right exact and
commutes with direct sums. The category of
$(\bf{G},\bf{F})$-bicomodules $\lrBicomod{\bf{G}}{\bf{F}}$ is a
preadditive category with cokernels and arbitrary direct sums.
\end{lemma}

Consider the categories of bicomodules $\lBicomod{\Bb}{\bf{F}}$ and
$\lrBicomod{\bf{G}}{\bf{F}}$. There are two functors connecting
those categories. The left forgetful functor $\sS:
\lrBicomod{\bf{G}}{\bf{F}} \to \lBicomod{\Bb}{\bf{F}}$, which sends
any $(\bf{G}, \bf{F})$-bicomodule $(\Sf{M},\mmm,\nnn)$ to the
$(\Bb,\bf{F})$-bicomodule $(\Sf{M},\mmm)$ and which is identical on
the morphisms. Secondly, the functor $\tT: \lBicomod{\Bb}{\bf{F}}
\to \lrBicomod{\bf{G}}{\bf{F}}$ which sends $(\Sf{M}',\mmm') \to
(\Sf{M'}G,\mmm'_G,\Sf{M}'\vartheta)$ and $\Sf{f} \to \Sf{f}_G$.
These functors form an adjunction, more precise we have

\begin{lemma}\label{adj}
For every pair of objects $\lr{(\Sf{N},\rrr,\sss),\,(M,\mmm)}$ of
$\lrBicomod{\bf{G}}{\bf{F}} \times \lBicomod{\Bb}{\bf{F}}$, there is
a natural transformation
$$\xymatrix@C=40pt@R=0pt{
\hom{\,\lrBicomod{\bf{G}}{\bf{F}}}{(\Sf{N},\rrr,\sss)}{\tT(\Sf{M},\mmm)}
\ar@{->}^-{\Phi_{\Sf{N},\,\Sf{M}}}[r] &
\hom{\lBicomod{\Bb}{\bf{F}}}{\sS(\Sf{N},\rrr,\sss)}{(\Sf{M},\mmm)}\\
\Sf{f} \ar@{|->}[r] & \Sf{M}\varsigma \circ \Sf{f} \\ \Sf{g}_G \circ
\sss & \Sf{g}. \ar@{|->}[l]}$$ That is $\sS$ is a left adjoint
functor to $\tT$.
\end{lemma}

Let $\Xx$ be the one-object category, then the category
$\lBicomod{\Xx}{\bf{F}}$ can be described as follows. A functor
$\Sf{X}:\Xx\to \Aa$ is completely determined by the image $X$ of the
single object in $\Xx$. A natural transformation $\xxx:\Sf{X}\to
F\Sf{X}$ is completely determined by a morphism $d^X:X\to F(X)$. In
this way, we can identify an object in $\lBicomod{\Xx}{\bf{F}}$ with
a pair $(X,d^X)$ consisting of an object $X \in \Aa$ and a morphism
$d^X: X \rightarrow F(X)$ satisfying $$ \delta_X \, \circ \, d^X
\,\, = \,\, F(d^X) \,\circ \, d^X,\quad \xi_X \, \circ \, d^X \,\,
=\,\,X.$$ Similarly, a morphism $f:(X,d^X) \rightarrow (X',d^{X'})$
in $\lBicomod{\Xx}{\bf{F}}$ is completely determined by a morphism
$f: X \rightarrow X'$ of $\Aa$ such that $$ d^{X'} \,\circ \, f \,\,
=\,\, F(f) \,\circ \, d^X.$$ Under this identification, we will
denote this category by $\Aa^{\bf F}$. Denote by ${\bf S}:\Aa^{\bf
F}\to \Aa$ the forgetful functor and ${\bf T}:\Aa\to \Aa^{\bf F},
{\bf T}(Y)=(F(Y),\delta_Y)$, ${\bf T}=F(f)$, for every object $Y$
and morphism $f$ of $\Aa$. Then we obtain an adjunction ${\bf S}
\dashv {\bf T}$, with ${\bf S}{\bf T} =F$ satisfying a universal
property, see \cite[Theorem 2.2]{Eilenberg/Moore:1965}.

\begin{rem}
It is well known that $\Aa^{{\bf F}}$ is an additive category with
direct sums and cokernels, admitting $(F(U),\delta_U)$ as a
sub-generator, whenever $U$ is a generator of $\Aa$. However,
$\Aa^{{\bf F}}$ is not necessarily a Grothendieck category. But, if
we assume that $F$ is an exact functor and that $\Aa$ poses a
generating set of finitely generated objects, then one can easily
check that $\Aa^{{\bf F}}$ becomes a Grothendieck category.
\end{rem}

The main results of this section is the following
\begin{theo}\label{Separable}
Let $\Aa$ be a Grothendieck category. Consider a comonad ${\bf
F}=(F,\delta,\xi)$ in $\Aa$ whose functor $F$ preserves cokernels
and commutes with direct sums. The following are equivalent
\begin{enumerate}[(i)]
\item ${\bf S}: \Aa^{{\bf F}} \longrightarrow \Aa$ is separable
functor;

\item $\sS:\, \lrBicomod{{\bf F}}{{\bf F}} \longrightarrow
\lBicomod{\Aa}{{\bf F}}$ is a Maschke functor;

\item $\delta: (F,\delta,\delta) \longrightarrow
(F^2,\delta_F,F\delta)$ is a split monomorphism in the category
$\lrBicomod{{\bf F}}{{\bf F}}$.
\end{enumerate}
\end{theo}
\begin{proof}
$(i) \Rightarrow (iii)$. The unit of the adjunction ${\bf S}
\dashv \bf{T}$ is given by
\begin{equation}\label{eta}
\xymatrix{\eta_{(X,\, d^X)}: (X,d^X) \ar@{->}^-{d^X}[rr] & & {\bf
TS}(X,d^X)\, =\, (F(X),\delta_X)}
\end{equation}
for every object $(X,d^X)$ of $\Aa^{\bf{F}}$. By hypothesis there is
a natural transformation $\psi: \bf{T} {\bf S} \rightarrow
\id_{\Aa^{\bf{F}}}$ such that $\psi \circ \eta \,=\,
\id_{\Aa^{\bf{F}}}$. Let us denote by $\nabla: F^2 \rightarrow F$
the natural transformation given by the collection of morphisms
$\nabla_X\,=\,{\bf S}(\psi_{(F(X), \,\delta_X)})$, where $X$ runs
through the class of object of $\Aa$. By construction $\nabla \circ
\delta =F$ and $\nabla: (F^2,\delta_F) \rightarrow (F,\delta)$ is a
morphism of the category $\lBicomod{\Aa}{{\bf F}}$. Since $\psi$ is
a natural transformation and $\delta_X: (F(X),\delta_X) \rightarrow
(F^2(X),\delta_{F(X)})$ is morphism in $\Aa^{\bf{F}}$, we have the
following commutative diagram $$ \xymatrix@C=60pt{ F^3
\ar@{->}^-{{\bf S}\psi_{F^2}}[r]  & F^2 \\ F^2 \ar@{->}^-{{\bf
S}\psi_F}[r] \ar@{->}^-{F\delta}[u]  & F\ar@{->}_-{\delta}[u] }$$
Therefore $\delta \circ \nabla\,=\, \nabla_F \circ F\delta$, which
means that $\nabla: (F^2,\delta_F,F\delta) \rightarrow
(F,\delta,\delta)$ is a morphism in the category $\lrBicomod{{\bf
F}}{{\bf F}}$. Thus $\delta$ is a split monomorphism of the category
$\lrBicomod{{\bf F}}{{\bf F}}$.
\\ $(iii) \Rightarrow (ii)$. Let us denote by $\Lambda:
(F^2,\delta_F,F\delta) \rightarrow (F,\delta,\delta)$ the left
inverse of $\delta: (F,\delta,\delta) \rightarrow
(F^2,\delta_F,F\delta)$, i.e. $\Lambda \circ \delta =F$, in the
category $\lrBicomod{{\bf F}}{{\bf F}}$. Let $(\Sf{M},\mmm,\nnn)$ be
any $\bf{F}$-bicomodule. The unit of the adjunction $\sS \dashv \tT$
stated in lemma \ref{adj}, at this bicomodule is given by
\begin{equation}\label{Theta}
\xymatrix@C=40pt{ \Theta_{(\Sf{M},\,\mmm,\,\nnn)}:
(\Sf{M},\mmm,\nnn) \ar@{->}^-{\nnn}[rr] & & \tT \circ
\sS(\Sf{M},\mmm,\nnn)\,=\, (MF,\mmm_F, M\delta).}
\end{equation} Consider the natural transformation defined by the following
composition
$$\xymatrix@C40pt{ \upsilon:\, \Sf{M}F \ar@{->}^-{\nnn_F}[r] &
\Sf{M}F^2 \ar@{->}^-{\Sf{M}\Lambda}[r] & \Sf{M}F
\ar@{->}^-{\Sf{M}\xi}[r] & \Sf{M} }.$$

It is easily seen that $ \upsilon \circ  \nnn \, =\, \Sf{M}$. The
implication will be established if we show that $\upsilon$ is a
morphism in the category of bicomodules $\lrBicomod{{\bf F}}{{\bf
F}}$. We can compute
\begin{eqnarray*}
  \mmm \circ \upsilon &=& \mmm \circ \Sf{M}\xi \circ \Sf{M}\Lambda \circ \nnn_F \\
   &=& F\Sf{M}\xi \circ \mmm_F \circ \Sf{M}\Lambda \circ \nnn_F,
   \qquad \,\, \mmm_{-}\,\, \text{ is } \,\,\text{ natural}  \\
   &=& F\Sf{M}\xi \circ F\Sf{M}\Lambda \circ \mmm_{F^2} \circ \nnn_F,
   \qquad \,\, \mmm_{-} \,\,\text{ is }\,\, \text{ natural} \\
   &=& F\Sf{M}\xi \circ F\Sf{M}\Lambda \circ F\nnn_F \circ \mmm_F, \qquad \,\,\text{ by}\,\,\eqref{comod-11} \\
   &=& F\left(\underset{}{} \Sf{M}\xi \circ \Sf{M}\Lambda \circ \nnn_F\right) \circ
   \mmm_F \\
   &=& F\upsilon \circ \nnn_F,
\end{eqnarray*}
which proves that $\upsilon$ is a morphism in
$\lBicomod{\Aa}{\bf{F}}$. On the other hand, we have
\begin{eqnarray*}
  \nnn \circ \upsilon  &=& \nnn \circ {\sf M}\xi \circ \Sf{M}\Lambda \circ \nnn_F \\
   &=& \Sf{M}F\xi \circ \nnn_F \circ \Sf{M}\Lambda \circ \nnn_F,
   \qquad \,\, \nnn_{-}\,\, \text{ is } \,\,\text{ natural} \\
   &=& \Sf{M}F\xi \circ \Sf{M}F\Lambda \circ \nnn_{F^2} \circ \nnn_F,
   \qquad \,\, \nnn_{-}\,\, \text{ is } \,\,\text{ natural} \\
  &=& \Sf{M}F\xi \circ \Sf{M}F\Lambda \circ \Sf{M}\delta_F \circ \nnn_F,
  \quad \,\, \text{ by } \,\,\eqref{comod-00} \\
   &=& \Sf{M}F\xi \circ \Sf{M}\delta \circ \Sf{M}\Lambda \circ  \nnn_F,
  \quad \,\, \text{ by }\,\, \eqref{bicolineal} \\
   &=&  \Sf{M}\Lambda \circ \nnn_F,
\end{eqnarray*} and
\begin{eqnarray*}
  \upsilon_F \circ \Sf{M}\delta &=& \Sf{M}\xi_F \circ \Sf{M}\Lambda_F  \circ \nnn_{F^2} \circ \Sf{M}\delta \\
   &=& \Sf{M}\xi \circ \Sf{M}\Lambda_F  \circ \Sf{M}F\delta \circ \nnn_F,
  \qquad \,\, \nnn_{-} \,\,\text{ is }\,\, \text{ natural} \\
   &=& \Sf{M}\xi_F \circ \Sf{M}(\Lambda_F  \circ F\delta) \circ \nnn_F \\
   &=& \Sf{M}\xi_F \circ \Sf{M}\delta \circ \Sf{M}\Lambda  \circ
   \nnn_F, \qquad \text{by}\,\, \eqref{bicolineal} \\
   &=& \Sf{M}\Lambda \circ \nnn_F.
\end{eqnarray*}
Therefore $\upsilon_F \circ \Sf{M}\delta\,=\,\nnn \circ \upsilon $
and $\upsilon$ is a morphism of $\bf{F}$-bicomodules. Hence
$\sS$ is a Maschke functor. \\
$(ii) \Rightarrow (i)$. Given $(\Sf{M},\mmm,\nnn)$ an
$\bf{F}$-bicomodule, we denote by
$$\xymatrix@C=40pt{\Gamma_{(\Sf{M},\,\mmm,\,\nnn)}: \tT
\sS(\Sf{M},\mmm,\nnn)\,=\, (\Sf{M}F,\mmm_F, \Sf{M}\delta)
\ar@{->}[r] & (\Sf{M},\mmm,\nnn)}$$ the splitting morphism of
$\Theta_{(\Sf{M},\mmm,\nnn)}$ in the category of
$\bf{F}$-bicomodules. Here $\Theta_{-}$ is the unit of the
adjunction $\sS \dashv \tT$. Since $(F,\delta,\delta)$ is
$\bf{F}$-bicomodule, we put
$\gamma:=\Gamma_{(F,\,\delta,\,\delta)}$, thus $\gamma \circ \delta
\,=\, F$. For any object $(X,d^X)$ of the category $\Aa^{\bf{F}}$,
we consider the composition
$$\xymatrix@C=40pt{\phi_{(X,\,d^X)}: F(X) \ar@{->}^-{F(d^X)}[r] &
F^2(X) \ar@{->}^-{\gamma_X}[r] & F(X) \ar@{->}^-{\xi_X}[r] & X}.$$
We claim that $\phi_{-}$ is a natural transformation which satisfies
$\phi_- \circ \eta_-\,=\, \id_{\Aa^{\bf{F}}}$, where $\eta_{-}$ is
the unit of the adjunction ${\bf S} \dashv \bf{T}$ given in
\eqref{eta}. First of all, we have
\begin{eqnarray*}
   \phi_{(X,\,d^X)} \circ \eta_{(X,\, d^X)} &=& \xi_X \circ \gamma_X
\circ F(d^X) \circ d^X \\
   &=& \xi_X \circ \gamma_X \circ  \delta_X \circ d^X \\
   &=& \xi_X \circ d^X \,\, =\,\, (X,d^X),
\end{eqnarray*}
for every object $(X,d^X)$ of $\Aa^{\bf{F}}$. To see that
$\phi_{(X,\,d^X)}$ is a morphism in $\Aa^{\bf{F}}$, we can compute
on one hand
\begin{eqnarray*}
  d^X \circ \phi_{(X,\,d^{X})} &=& d^X \circ \xi_X \circ \gamma_X \circ F(d^X) \\
   &=& \xi_{F(X)} \circ F(d^X) \circ \gamma_X \circ F(d^X), \qquad \xi_{-}\,\, \text{is} \,\, \text{natural} \\
   &=& \xi_{F(X)} \circ \gamma_{F(X)} \circ F^2(d^X) \circ F(d^X),\qquad \gamma_{-}\,\, \text{is} \,\, \text{natural} \\
   &=& \xi_{F(X)} \circ \gamma_{F(X)} \circ F\delta \circ F(d^X) \\
   &=& \xi_{F(X)} \circ \delta_X \circ \gamma_{X} \circ F(d^X),\qquad
   \text{by}\,\,\eqref{bicolineal}\\
   &=& \gamma_X \circ F(d^X)
\end{eqnarray*} and secondly,
\begin{eqnarray*}
  F\phi_{(X,\, d^X)} \circ \delta_X &=& F\xi_X \circ F\gamma_X \circ F^2(d^X) \circ \delta_X \\
   &=& F\xi_X \circ F\gamma_X \circ \delta_{F(X)} \circ F(d^X),\qquad \delta_{-}\,\, \text{is} \,\, \text{natural}  \\
   &=& F\xi_X \circ \delta_X \circ \gamma_{X} \circ F(d^X),\qquad
   \text{by}\,\,\eqref{bicolineal} \\
   &=& \gamma_X \circ F(d^X).
\end{eqnarray*}
Therefore, $F\phi_{(X,\,d^X)} \circ \delta_X\,=\,d^X \circ
\phi_{(X,\,d^{X})}$. Lastly, if we consider a morphism $f:(X,d^X)
\rightarrow (Y,d^Y)$ in $\Aa^{\bf{F}}$, then
\begin{eqnarray*}
  f \circ \phi_{(X,\, d^X)} &=& f \circ \xi_X \circ \gamma_X \circ F(d^X)  \\
   &=& \xi_Y \circ F(f) \circ \gamma_X \circ F(d^X),\qquad \xi_{-}\,\, \text{is} \,\, \text{natural}  \\
   &=& \xi_Y \circ \gamma_Y \circ F^2(f) \circ F(d^X), \qquad \gamma_{-}\,\, \text{is} \,\, \text{natural} \\
   &=& \xi_Y \circ \gamma_Y \circ F(d^Y) \circ F(f) \\
   &=& \phi_{(Y,\,d^Y)} \circ F(f),
\end{eqnarray*} which shows that $\phi_{-}$ is a natural
transformation.
\end{proof}

\section{Coderivations and Cointegrations}\label{Sect-2}

Let ${\bf F}=(F,\delta,\xi)$ be a comonad in $\Aa$ with underlying
functor $F \in \ContFunt{\Aa}{\Aa}$. Consider a bicomodule
$(\Sf{M},\mmm,\nnn) \in \lrBicomod{{\bf F}}{{\bf F}}$. A
\emph{coderivation} from $\Sf{M}$ to $ F$ is a natural
transformation $ \Sf{g}: \Sf{M} \longrightarrow F$ such that
\begin{equation}\label{Coderv}
\delta  \circ \Sf{g} \,\, = \,\, F\Sf{g} \circ \mmm \, +\,
\Sf{g}_{F} \circ \nnn.
\end{equation}

The set of all coderivations from $(\Sf{M},\mmm,\nnn)$ is an
additive group which we denote by $\Coder{\Sf{M}}{F}$. A
coderivation $g \in \Coder{\Sf{M}}{F}$ is said to be \emph{inner}
if there exists a natural transformation $\lambda: \Sf{M}
\rightarrow \id_{\Aa}$ such that
\begin{equation}\label{coder-int}
g \,\,=\,\, \lambda_F \circ \nnn \,-\, F\lambda \circ \mmm.
\end{equation} The sub-group of all inner coderivations will be
denoted by $\IntCoder{\Sf{M}}{F}$.

Let $(\Sf{M},\mmm,\nnn)$ and $(\Sf{M}',\mmm',\nnn')$ be two ${\bf
F}$-bicomodules. A \emph{left cointegration} from
$(\Sf{M},\mmm,\nnn)$ into $(\Sf{M}',\mmm',\nnn')$ is a natural
transformation ${\sf h}: {\sf M} \rightarrow {\sf M}'F$ which
satisfies
\begin{equation}\label{Cointegracion}
\mmm'_F \circ \Sf{h} \,\,=\,\, F \Sf{h} \circ \mmm,\quad
\Sf{M}'\delta \circ \Sf{h} \,\, =\,\, \nnn'_F \circ \Sf{h} \,+\,
\Sf{h}_{F}\circ \nnn.
\end{equation} The first equality means that ${\sf h}: \sS(\Sf{M},\mmm,\nnn)=({\sf M},\mmm)
\rightarrow \sS\tT\sS(\Sf{M},\mmm,\nnn)=({\sf M}'F,\mmm'_F)$ is a
morphism in the category $\lBicomod{\Aa}{\bf{F}}$. Right
cointegrations are defined in a similar way. Since we are only
concerned with the left ones, we will not mention the word
``left'' before cointegration. The additive group of all
cointegrations from $(\Sf{M},\mmm,\nnn)$ into
$(\Sf{M}',\mmm',\nnn')$ will be denoted by
$\Coint{\Sf{M}}{\Sf{M}'}$. A cointegration ${\sf h} \in
\Coint{\Sf{M}}{\Sf{M}'}$ is said to be \emph{inner} if there
exists a natural transformation $\varphi: \Sf{M} \rightarrow
\Sf{M}'$ which satisfies
\begin{equation}\label{Coint-Int}
\mmm' \circ \varphi \,\, =\,\, F\varphi \circ \mmm, \quad
\Sf{h}\,\,=\,\, \varphi_F \circ \nnn\,-\, \nnn' \circ \varphi.
\end{equation} The first equality means that $\varphi:(\Sf{M}, \mmm)
\rightarrow (\Sf{M}',\mmm')$ is a morphism in the category
$\lBicomod{\Aa}{\bf{F}}$. The sub-group of all inner
cointegrations will be denoted by $\IntCoint{\Sf{M}}{\Sf{M}'}$.
The following proposition was first stated for bimodule over ring
extension in \cite{Kleiner:1985} and for bicomodules over corings
in \cite{Guzman:1989}. For the sake of completeness, we give the
proof.
\begin{prop}\label{Coint-Coder}
For $(\Sf{M}, \mmm,\nnn)$ any $\bf{F}$-bicomodule, there is a
natural isomorphism of additive groups
$$\xymatrix@C=60pt@R=0pt{
\Coint{\Sf{M}}{F} \ar@{->}^-{\sim}[r] & \Coder{\Sf{M}}{F} \\
{\sf h} \ar@{|->}[r] & \xi_F \circ{\sf h} \\ F{{\sf g}} \circ\mmm
& \ar@{|->}[l] {\sf g}}$$ whose restriction to the inner
sub-groups gives again an isomorphism $$\IntCoint{\Sf{M}}{F} \cong
\IntCoder{\Sf{M}}{F}.$$
\end{prop}
\begin{proof}
We only show that the mutually inverse maps are well defined. Let
$\Sf{h} \in \Coint{\Sf{M}}{F}$, and put $\Sf{g}:= \xi_F \circ
\Sf{h}$. We have
\begin{eqnarray*}
  \delta \circ g &=& \delta \circ \xi_F \circ \Sf{h} \\
   &=& \xi_{F^2} \circ F\delta \circ \Sf{h}, \qquad
   \delta_{-}\,\, \text{is} \,\, \text{natural} \\
   &=& \xi_{F^2} \circ \left(\underset{}{} \delta_F \circ \Sf{h} + \Sf{h}_F \circ \nnn \right) \\
   &=& \left(\xi_F \circ \delta\right)_F \circ \Sf{h} \,+\, \xi_{F^2} \circ \Sf{h}_F \circ \nnn \\
   &=& \Sf{h} \,+\,\xi_{F^2} \circ \Sf{h}_F \circ \nnn
\end{eqnarray*} and
$$F\xi_F \circ F\Sf{h} \circ \mmm \,+\, \xi_{F^2} \circ \Sf{h}_F \circ
\nnn\,\,=\,\,F \xi_F \circ \delta_F \circ \Sf{h} \,+\, \xi_{F^2}
\circ \Sf{h}_F \circ\nnn \,\,=\,\, \Sf{h} + \xi_{F^2} \circ
\Sf{h}_F \circ \nnn.$$ That is $\Sf{g} \in \Coder{\Sf{M}}{F}$.
Conversely, given $\Sf{g} \in \Coder{\Sf{M}}{F}$, we put
$\Sf{h}\,=\, F\Sf{g} \circ \mmm$. We find
\begin{eqnarray*}
  \delta_F \circ \Sf{h} &=& \delta_F \circ F\Sf{g} \circ \mmm \\
   &=& {F^2}\Sf{g} \circ \delta_{\Sf{M}} \circ \mmm, \qquad \delta_{-} \,\, \text{is}\,\, \text{natural} \\
   &=& {F^2}\Sf{g} \circ F\mmm \circ \mmm, \qquad \text{by} \,\, \eqref{comod-00} \\
   &=& F \Sf{h} \circ \mmm,
\end{eqnarray*} which shows the first equality of equation
\eqref{Cointegracion}. Now,
\begin{eqnarray*}
  F\delta \circ \Sf{h} &=& F\delta \circ F\Sf{g} \circ \mmm \\
   &=& F\delta \circ \delta \circ  \Sf{g} - F\delta \circ \Sf{g}_{F} \circ \nnn  \\
   &=& \delta_F \circ \delta \circ \Sf{g} - \Sf{g}_{F^2} \circ \Sf{M}\delta \circ\nnn, \qquad
   \Sf{g}_{-} \,\, \text{is}\,\, \text{natural} \\
   &=& \delta_F\circ \delta \circ \Sf{g} -  \Sf{g}_{F^2} \circ \nnn_F \circ \nnn,
   \qquad \text{by}\,\, \eqref{comod-11} \,\, \text{and}\,\, \eqref{comod-00} \\
   &=& \delta_F\circ \delta \circ \Sf{g}-\delta_F \circ \Sf{g}_F \circ \nnn +
   \delta_F \circ \Sf{g}_F \circ \nnn -  \Sf{g}_{F^2} \circ \nnn_F \circ \nnn   \\
   &=& \delta_F\circ \left(\underset{}{}\delta \circ \Sf{g}-  \Sf{g}_F \circ \nnn\right) +
   \left(\underset{}{} \delta_F \circ \Sf{g}_F  -  \Sf{g}_{F^2} \circ \nnn_F\right) \circ \nnn \\
   &=& \delta_F\circ \left(\underset{}{}\delta \circ \Sf{g}-  \Sf{g}_F \circ \nnn\right) +
   \left(\underset{}{} \delta \circ \Sf{g}  -  \Sf{g}_{F} \circ \nnn\right)_F \circ \nnn \\
   &=& \delta_F\circ F \Sf{g} \circ \mmm + F\Sf{g}_F \circ \mmm_F \circ
   \nnn \\
   &=& \delta_F \circ \Sf{h} + \Sf{h}_F \circ \nnn
\end{eqnarray*} which proves that $\Sf{h}=F\Sf{g} \circ \mmm \in
\Coint{\Sf{M}}{F}$.
\end{proof}

Following \cite{Guzman:1989}, we will give in the next step the
notion of universal cointegration and that of universal
coderivation.

Given $(\Sf{M},\mmm,\nnn)$ any $\bf{F}$-bicomodule, consider the
$\bf{F}$-bicomodule $({\sf M}F,\mmm_F,\Sf{M}\delta)$, which is the
image of $(\Sf{M},\mmm,\nnn)$ under the functor $\tT \sS$. We call
it the bicomodule induced by $\Sf{M}$. Since $\nnn: ({\sf
M},\mmm,\nnn) \rightarrow ({\sf M}F,\mmm_F,{\sf M}\delta)$ is a
morphism of $\bf{F}$-bicomodules, we obtain by lemma \ref{conucleos}
the following sequence of $\bf{F}$-bicomodules
\begin{equation}\label{sucesion}
\xymatrix@C=20pt{ 0 \ar@{->}[r] & ({\sf M},\mmm,\nnn)
\ar@{->}^-{\nnn}[r] & ({\sf M}F,\mmm_F,{\sf M}\delta)
\ar@{->}^-{\nnn^c}[r] & (\kK({\sf M}), \uuu,\vvv) \ar@{->}[r] & 0},
\end{equation}
where $(\kK({\sf M}), \uuu,\vvv),\nnn^c$ denotes the cokernel of
$\nnn$ in the category $\lrBicomod{\bf F}{\bf F}$. Notice , that
this still be a cokernel in the category $\lBicomod{\Aa}{\bf F}$,
after forgetting by $\sS$. Consider the natural transformation
$${\sf w}':= {\sf M}F \,-\, \nnn \circ {\Sf M}\xi\,: {\sf M}F
\longrightarrow {\sf M}F$$ It is easily checked that $\mmm_F \circ
{\sf w}'\,=\,F{\sf w}' \circ \mmm_F$, thus ${\sf w}'$ is a
morphism in the category $\lBicomod{\Aa}{\bf{F}}$. Also, $\Sf{w}'$
satisfies $\Sf{w}' \circ \nnn \,=\,0$. So, by the universal
property of cokernels, there exists a morphism in the category
$\lBicomod{\Aa}{\bf{F}}$, $\Sf{w}: (\kK(\Sf{M}), \vvv) \rightarrow
(\Sf{M}F,\mmm_F)$ which makes the following diagram commutative
\begin{equation}\label{omegas}
\xymatrix@C=50pt{ \Sf{M} \ar@{->}^-{\nnn}[r] & \Sf{M}F
\ar@{->}^-{\nnn^c}[r] \ar@{->}_-{\Sf{w}'}[d] & \kK(\Sf{M}) \ar@{-->}^-{\Sf{w}}[dl]\\
& \Sf{M}F & }
\end{equation}
Thus $\Sf{w} \circ \nnn^c \,=\, \Sf{w}'$, and so $\nnn^c \circ
\Sf{w} \circ \nnn^c\,=\, \nnn^c$. Hence $\nnn^c \circ \Sf{w} \,=\,
\kK(\Sf{M})$, since $\nnn^c$ is an epimorphism. Furthermore, we
have
\begin{prop}\label{univ-coint}
The morphism $\Sf{w}$ is a cointegration into $\Sf{M}$ (i.e. $\Sf{w}
\in \Coint{\kK(\Sf{M})}{\Sf{M}}$) which satisfies the following
universal property. For every $\bf{F}$-bicomodule $(\Sf{M}',\mmm'
,\nnn')$ and every cointegration $\Sf{h} \in
\Coint{\Sf{M}'}{\Sf{M}}$, there exists a morphism of
$\bf{F}$-bicomodules $\Sf{f}: (\Sf{M}',\mmm',\nnn') \rightarrow
(\kK(\Sf{M}), \uuu,\vvv)$ such that $ \Sf{h}\,\,=\,\ \Sf{w} \circ
\Sf{f}$. Moreover, the following are equivalent
\begin{enumerate}[(i)]
\item The sequence $$\xymatrix@C=30pt{ 0 \ar@{->}[r] &
(\Sf{M},\mmm,\nnn) \ar@{->}^-{\nnn}[r] &
(\Sf{M}F,\mmm_F,\Sf{M}\delta) \ar@{-}^-{\nnn^c}[r] &
(\kK(\Sf{M}),\uuu,\vvv) \ar@{->}[r] & 0 }$$ splits in the category
of bicomodules $\lrBicomod{\bf{F}}{\bf{F}}$.

\item The universal cointegration $\Sf{w}: \kK(\Sf{M}) \to \Sf{M}F$
is inner.
\end{enumerate}
\end{prop}
\begin{proof}
For the first statement, it is enough to show that $\Sf{w}'$ is
cointegration into $\Sf{M}$, since $\nnn^c$ is an epimorphism. By
definition $\Sf{w}'$ satisfies the first equality in
\eqref{Cointegracion}. The second equality in \eqref{Cointegracion},
is given as follows  $$\Sf{M}\delta \circ \Sf{w}' \,\,=\,\,
\Sf{M}\delta - \Sf{M}\delta \circ \nnn \circ \Sf{M}\xi \,\,=\,\,
\Sf{M}\delta -\nnn_F \circ \nnn \circ \Sf{M}\xi$$ and
\begin{eqnarray*}
  \nnn_F \circ \Sf{w}' + \Sf{w}'_F \circ \Sf{M}\delta &=& \nnn_F - \nnn_F \circ \nnn \circ \Sf{M}\xi +
\Sf{M}\delta - \nnn_F \circ \Sf{M}\xi_F \circ \Sf{M}\delta \\
   &=& {\Sf{M}}\delta  -\nnn_F\circ \nnn \circ \Sf{M}\xi  \\
   &=& {\Sf{M}}\delta \circ \Sf{w}'
\end{eqnarray*}
The fact that $\Sf{w}$ is universal follows from the following
isomorphism of additive groups
\begin{equation}\label{Isomofos}
\xymatrix@R=0pt@C=40pt{
\hom{\,\,\lrBicomod{\bf{F}}{\bf{F}}}{\Sf{M}'}{\kK(\Sf{M})}
\ar@{->}^-{\sim}[r] & \Coint{\Sf{M}'}{\Sf{M}} \\ \varphi
\ar@{|->}[r] & \Sf{w} \circ \varphi  \\ \nnn^c \circ \Sf{h}  &
\Sf{h} \ar@{|->}[l] }
\end{equation}
whose proof is an easy computation. Now we check the equivalent statements. \\
$(i) \Rightarrow (ii)$. Let us denote by $\lambda:
(\kK(\Sf{M}),\uuu,\vvv) \to (\Sf{M}F,\mmm_F,\Sf{M}\delta)$ the right
inverse of $\nnn^c$ in the category $\lrBicomod{\bf{F}}{\bf{F}}$,
i.e. $\nnn^c \circ \lambda\,=\, \kK(\Sf{M})$. Define the composition
$$\xymatrix@C=30pt{\varphi:\,\, \kK(\Sf{M}) \ar@{->}^-{\lambda}[r] &
\Sf{M}F \ar@{->}^-{\Sf{M}\xi}[r] & \Sf{M} }.$$ Then we have
$$\mmm\circ \varphi\,=\, \mmm \circ\Sf{M}\xi \circ \lambda\,=\,
F\Sf{M}\xi \circ \mmm \circ \lambda \,\,=\,\, F\Sf{M}\xi \circ
F\lambda \circ \uuu\,\,=\,\, F\left(\underset{}{}\Sf{M}\xi \circ
\lambda \right) \circ \uuu \,\,=\,\, F\varphi \circ \uuu,$$ which
entails that $\varphi$ is morphism in $\lBicomod{\Aa}{\bf{F}}$. The
cointegration $\Sf{w}$ is inner by $\varphi$. Namely,
\begin{eqnarray*}
  \varphi_F \circ \vvv - \nnn \circ \varphi &=& \Sf{M}\xi_F
  \circ \lambda_F \circ \vvv - \nnn \circ \Sf{M}\xi \circ \lambda \\
   &=& \Sf{M}\xi_F \circ \Sf{M}\delta \circ \lambda - \nnn \circ \Sf{M}\xi \circ \lambda \\
   &=& \lambda - \nnn \circ \Sf{M} \xi  \circ \lambda  \\
   &=& \left(\underset{}{} \Sf{M}F - \nnn \circ \Sf{M}\xi \right) \circ \lambda  \\
   &=& \Sf{w} \circ \nnn^c \circ \lambda \,\, =\,\, \Sf{w}.
\end{eqnarray*}
$(ii) \Rightarrow (i)$. Suppose that there exists $\beta:
\kK(\Sf{M}) \to \Sf{M}$ a morphism in $\lBicomod{\Aa}{\bf{F}}$ such
that $\Sf{w}\,=\, \beta_F \circ\vvv - \nnn \circ \beta$. Consider
the natural transformation
$$\xymatrix@C=30pt{ \Gamma:\, \kK(\Sf{M}) \ar@{->}^-{\vvv}[r]&
\kK(\Sf{M})F \ar@{->}^-{\beta_F}[r] & \Sf{M}F}.$$ Then we find
$\nnn^c \circ \Gamma \,\,=\,\, \nnn^c \circ \beta_F \circ \vvv
\,\,=\,\, \nnn^c \circ \Sf{w} + \nnn^c \circ \nnn \circ \beta
\,\,=\,\, \nnn^c \circ \Sf{w}\,\,=\,\, \kK(\Sf{M})$. Furthermore,
$\Gamma$ is a morphism in the category of bicomodules
$\lrBicomod{\bf{F}}{\bf{F}}$, as the following commutative diagrams
shown
$$\xymatrix@C=25pt{\kK(\Sf{M}) \ar@{->}^-{\vvv}[r] \ar@{->}_-{\uuu}[d] &
\kK(\Sf{M})F \ar@{->}^-{\beta_F}[r] \ar@{->}_-{\uuu_F}[d] & \Sf{M}F
\ar@{->}^-{\mmm_F}[d] \\ F\kK(F) \ar@{->}_-{F\vvv}[r] &
F\kK(\Sf{M})F \ar@{->}_-{F\beta_F}[r]& F\Sf{M}F} \qquad
\xymatrix@C=40pt{\kK(\Sf{M}) \ar@{->}^-{\vvv}[r] \ar@{->}_-{\vvv}[d]
&\kK(\Sf{M})F \ar@{->}^-{\beta_F}[r]
\ar@{->}_-{\kK(\Sf{M})\delta}[d]
& \kK(\Sf{M})F \ar@{->}^-{\Sf{M}\delta}[d] \\
\kK(\Sf{M})F \ar@{->}_-{\vvv_F}[r] & \kK(\Sf{M})F^2
\ar@{->}_-{\beta_{F^2}}[r]& \Sf{M}F^2}$$ Therefore the sated
sequence splits in the category $\lrBicomod{\bf{F}}{\bf{F}}$.
\end{proof}

The cointegration $\Sf{w}$ from Proposition \ref{univ-coint} will
be referred to as the \emph{universal cointegration into}
$\Sf{M}$.

From now on $\Sf{w}$ denotes the universal cointegration into the
$\bf{F}$-bicomodule $(F,\delta,\delta)$. That is $\Sf{w}: \kK(F)
\rightarrow F^2$ with properties $\Sf{w} \circ \delta^c \,\,=\,\,
F^2 - \delta \circ F\xi$ and $\delta^c \circ \Sf{w} \,=\, \kK(F)$,
where
$$\xymatrix@C=20pt{0 \ar@{->}[r] & (F,\delta,\delta)
\ar@{->}^-{\delta}[r]& (F^2,\delta_F,F\delta)
\ar@{->}^-{\delta^c}[r] & (\kK(F),\uuu,\vvv) \ar@{->}[r] &0}$$ is
the canonical sequence. Consider the natural transformation
$\Sf{d}: \kK(F) \rightarrow F$ defined by $\Sf{d}: \,\,=\,\, F\xi
\circ \Sf{w}\,-\, \xi_F \circ \Sf{w}$.

\begin{lemma}\label{uni-coder}
The morphism $\Sf{d}$ is a coderivation with the following universal
property. For every $\bf{F}$-bicomodule $(\Sf{M},\mmm,\nnn)$ and
every coderivation $\Sf{g} \in \Coder{\Sf{M}}{F}$, there exists a
natural transformation $g': \Sf{M} \rightarrow \kK(F)$ such that
$\Sf{d} \circ g' \,=\, \Sf{g}$.
\end{lemma}
\begin{proof}
On one hand, we have \begin{eqnarray*}
  \delta \circ \Sf{d} \circ \delta^c &=& \delta \circ F\xi \circ \Sf{w} \circ \delta^c -
  \delta \circ \xi_F \circ \Sf{w} \circ \delta^c \\
   &=& \delta \circ F\xi - \delta \circ F\xi \circ \delta \circ F\xi -
  \delta \circ \xi_F + \delta \circ \xi_F \circ \delta \circ F\xi \\
   &=& \delta \circ F\xi - \delta \circ  F\xi - \delta \circ \xi_F + \delta \circ F\xi \\
   &=& -\delta \circ \xi_F + \delta \circ F\xi
\end{eqnarray*} on the other hand, we have
\begin{eqnarray*}
  \left(\underset{}{} F\Sf{d} \circ \uuu + \Sf{d}_F \circ \vvv \right) \circ \delta^c &=&
  F\Sf{d} \circ \uuu \circ \delta^c + \Sf{d}_F \circ \vvv \circ \delta^c  \\
   &=& F^2\xi \circ F\Sf{w} \circ \uuu \circ \delta^c - F\xi_{F} \circ F\Sf{w} \circ\uuu \circ \delta^c \\
   &\,\,&  + F\xi_F \circ \Sf{w}_F \circ\vvv\circ \delta^c  - \xi_{F^2}\circ \Sf{w}_F \circ\vvv\circ \delta^c \\
   &=& F^2\xi \circ F\Sf{w} \circ F\delta^c \circ \delta_F - F\xi_F \circ  F\Sf{w} \circ F\delta^c \circ
   \delta_F \\ &\,\,&
    + F\xi_{F} \circ \Sf{w}_F \circ \delta^c_F \circ F\delta -
    \xi_{F^2} \circ \Sf{w}_F \circ \delta^c_F \circ F\delta \\
   &=&  F^2\xi \circ F\left(\underset{}{} \Sf{w} \circ \delta^c \right) \circ
   \delta_F - F\xi_F \circ  F\left(\underset{}{}\Sf{w} \circ \delta^c\right) \circ
   \delta_F \\ &\,\,& + F\xi_{F} \circ \left(\underset{}{} \Sf{w} \circ \delta^c\right)_F \circ F\delta -
   \xi_{F^2} \circ \left(\underset{}{}\Sf{w} \circ \delta^c\right)_F \circ F\delta \\
   &=& F^2\xi \circ F\left(\underset{}{} F^2 - \delta \circ F\xi \right) \circ
   \delta_F - F\xi_F \circ  F\left(\underset{}{}F^2 - \delta \circ F\xi\right) \circ
   \delta_F \\ &\,\,& + F\xi_{F} \circ \left(\underset{}{} F^2 - \delta \circ F\xi\right)_F \circ F\delta -
   \xi_{F^2} \circ \left(\underset{}{}F^2 - \delta \circ F\xi\right)_F \circ F\delta \\
  &=& F^2 \xi \circ \delta_F - \xi_{F^2} \circ F \delta \,\,=\,\,
  \delta\circ F\xi - \delta \circ \xi_F \,\,=\,\, \delta \circ
  \Sf{d} \circ \delta^c,
\end{eqnarray*}
thus $\delta \circ \Sf{g}\,=\,F\Sf{d} \circ \uuu + \Sf{d}_F \circ
\vvv$, which shows that $\Sf{d} \in \Coder{\kK(F)}{F}$. Let now
$\Sf{g} \in \Coder{\Sf{M}}{F}$ be any coderivation. We know by
proposition \ref{Coint-Coder}, that $F\Sf{g} \circ \mmm \in
\Coint{\Sf{M}}{F}$. Using the isomorphism stated in
\eqref{Isomofos}, we obtain the following equality $$\Sf{w} \circ
\delta^c \circ F\Sf{g}\circ \mmm \,\, =\,\, F\Sf{g} \circ \mmm,$$
which implies that
$$\Sf{g}\,\, =\,\, \xi_F \circ \Sf{w} \circ \delta^c \circ F\Sf{g}
\circ \mmm,$$ as $\Sf{g}_{-}$ is natural. Developing $\Sf{d} \circ
\delta^c \circ \Sf{g}_F \circ \nnn$, we get
\begin{eqnarray*}
  \Sf{d} \circ \delta^c \circ \Sf{g}_F \circ \nnn &=&
  \left(\underset{}{} F\xi - \xi_F\right) \circ \Sf{w} \circ \delta^c \circ \Sf{g}_F \circ
  \nnn \\
   &=& \left(\underset{}{} F\xi - \xi_F\right) \circ \left(\underset{}{} F^2 - \delta \circ \xi_F\right) \circ
   \Sf{g}_F \circ \nnn \\
   &=& \left( \underset{}{} F\xi- F\xi \circ \delta \circ F\xi - \xi_F + \xi_F \circ \delta \circ F\xi\right) \circ
   \Sf{g}_F \circ \nnn \\
   &=& \left( \underset{}{} F\xi -F\xi -\xi_F + F\xi \right) \circ \Sf{g}_F \circ \nnn  \\
   &=& \left( \underset{}{} F\xi -\xi_F \right) \circ \Sf{g}_F \circ \nnn  \\
   &=& \left( \underset{}{} F\xi -\xi_F \right) \circ
   \left(\underset{}{} \delta \circ \Sf{g} - F\Sf{g} \circ\mmm \right) \\
   &=& -F\xi \circ F\Sf{g} \circ \mmm + \xi_F\circ F\Sf{g} \circ \mmm \\
   &=& \xi_F \circ \left( \underset{}{} F^2 - \delta \circ F\xi\right) \circ F\Sf{g} \circ
   \mmm \\
   &=& \xi_F \circ \Sf{w} \circ \delta^c \circ F\Sf{g} \circ
   \mmm  \,\, =\,\, \Sf{g}.
\end{eqnarray*} If we take $g'=\delta^c \circ F\Sf{g} \circ \mmm$,
then we find $\Sf{g}\,=\, \Sf{d} \circ g'$ and the universal
property is fulfilled.
\end{proof}

\begin{cor}\label{split-sequence}
Let $\bf{F}=(F,\delta,\delta)$ be a comonad in a Grothendieck
category $\Aa$ such that $F$ is a right exact and commutes with
direct sums. Consider the universal cointegration $\Sf{w}$ and the
universal coderivation $\Sf{d}$ associated to the
$\bf{F}$-bicomodule $(F,\delta,\delta)$. The following are
equivalent
\begin{enumerate}[(i)]
\item The sequence $$\xymatrix@C=20pt{0 \ar@{->}[r] & F
\ar@{->}^-{\delta}[r]& F^2 \ar@{->}^-{\delta^c}[r] & \kK(F)
\ar@{->}[r] &0}$$ is a split sequence in the category of
bicomodules $\lrBicomod{\bf{F}}{\bf{F}}$.

\item The universal cointegration $\Sf{w}$ is inner.

\item The universal coderivation $\Sf{d}$ is inner.
\end{enumerate}
\end{cor}
\begin{proof} The equivalence $(i) \Leftrightarrow (ii)$ is consequence of Proposition
\ref{univ-coint}. Let us check the equivalence between $(ii)$ and
$(iii)$.\\
$(ii) \Rightarrow (iii)$. We know there exists a morphism $\varphi:
\kK(F) \to F$ in $\lBicomod{\Aa}{\bf F}$ such that
$\Sf{w}\,=\, \varphi_F \circ \vvv - \delta \circ \varphi$. We
have
\begin{eqnarray*}
  \xi_F \circ \Sf{w} &=& \xi_F \circ \varphi_F \circ \vvv - \varphi
  \,\,=\,\, \lr{\xi \circ \varphi}_F \circ \vvv - \varphi \\
  F\xi \circ \Sf{w} &=& F\xi \circ \varphi_F \circ \vvv - \varphi.
\end{eqnarray*}
Hence
\begin{eqnarray*}
  \Sf{d} &=& F\xi \circ \varphi_F \circ \vvv - \varphi - \xi_F \circ \varphi_F \circ \vvv + \varphi\\
   &=& \varphi \circ \kK(F)\xi \circ \vvv - \xi_F \circ \varphi_F \circ \vvv \\
   &=& \varphi - \xi_F \circ \varphi_F \circ \vvv.
\end{eqnarray*}
But $F\lr{\xi \circ \varphi} \circ \uuu \,=\,
\varphi$, as $\varphi$ is a morphism in $\lBicomod{\Aa}{\bf F}$,
which proves that $\Sf{d}$ is inner by $-\xi \circ \varphi$.\\
$(iii) \Rightarrow (ii)$. Let us denote by $\lambda:\kK(F) \to
\id_{\Aa}$ the natural transformation which satisfies $\Sf{d}\,=\,
\lambda_F \circ \vvv - F\lambda \circ \uuu$. We define the map $\psi$ as the following composition
$\psi\,=\, F\lambda \circ \uuu: \kK(F) \to F\kK(F) \to F$.
This $\psi$ satisfies
$$ \delta \circ \psi \,\,=\,\, \delta \circ F\lambda \circ \uuu
\,\,=\,\, F^2\lambda \circ \delta_{\kK(F)} \circ \uuu \,\,=\,\,
F^2\lambda \circ F\uuu \circ \uuu \,\, =\,\, F\psi \circ   \uuu,$$
that is $\psi$ is a morphism in $\lBicomod{\Aa}{\bf F}$. The
universal cointegration is inner by $-\psi$, as the following
computations show
\begin{eqnarray*}
  \psi_F \circ \vvv - \delta \circ \psi &=& \lr{F\lambda \circ \uuu}_F \circ \vvv - \delta \circ F\lambda \circ \uuu \\
   &=& F\lambda_F \circ \uuu_F \circ \vvv - \delta \circ F\lambda \circ \uuu \\
   &=& F\lambda_F \circ F\vvv \circ \uuu - \delta \circ F\lambda \circ \uuu  \\
   &=& F\lr{\lambda_F \circ \vvv} \circ \uuu - \delta \circ F\lambda \circ \uuu \\
   &=& F\lr{\lr{F\xi - \xi_F} \circ \Sf{w} + F\lambda \circ \uuu} \circ \uuu - \delta \circ F\lambda \circ \uuu \\
   &=& F^2\xi \circ F\Sf{w} \circ \uuu - F\xi_F \circ F\Sf{w} \circ \uuu + F^2 \lambda \circ F\uuu \circ \uuu
   - \delta \circ F\lambda \circ \uuu \\
   &=& F^2\xi \circ \delta_F \circ \Sf{w} - F\xi_F \circ \delta_F \circ \Sf{w} + \delta \circ F\lambda \circ \uuu
   - \delta \circ F\lambda \circ \uuu \\
   &=& F^2\xi \circ \delta_F \circ \Sf{w} - \lr{F\xi \circ \delta}_F
   \circ
   \Sf{w}\\
   &=& \delta \circ F\xi \circ \Sf{w} -\Sf{w} \,\, =\,\,
   \lr{\delta \circ F\xi - F^2} \circ \Sf{w}\,\,=\,\, -\Sf{w}\circ
   \delta^c \circ \Sf{w}\,\,=\,\, -\Sf{w}.
\end{eqnarray*}

\end{proof}

\section{Cohomology For Bicomodules}\label{Sect-3}

The following lemma which will be used in the sequel, was in part
proved in \cite[Theorem 3.4]{Caenepeel/Militaru:2003}.
\begin{lemma}\label{lema-Sinj}
Let $\aA$ and $\bB$ two preadditive categories with cokernels, and
$\fF: \aA \to \bB$ a covariant functor with right adjoint functor
$\gG: \bB \to \aA$. Denote by $\chi$ and $\theta$ respectively,
the counit and unit of this adjunction. Let $\eE_0$ be the
injective class of all cosplit sequences in $\bB$, and put
$\eE=\fF^{-1}(\eE_0)$. For every object $M \in \aA$, the following
are equivalent
\begin{enumerate}[(i)]
\item $M$ is $\fF$-injective.

\item $M$ is $\eE$-injective.

\item The unit at $M$, $\theta_M: M \to \gG\fF(M)$, is a
split-mono in $\aA$.
\end{enumerate}
In particular every object of the form $\gG(N)$ is
$\eE$-injective, for every object $N \in \bB$. The functor $\fF$
is then Maschke if and only if the class of $\eE$-injective
objects coincides with class of all objects of $\aA$.
\end{lemma}
\begin{proof}
$(i) \Rightarrow (iii)$. We known by adjunction properties that
$\chi_{\fF(M)}\circ \fF(\theta_M)\,=\, \fF(M)$. Since
$M$ is $\fF$-injective, $\theta_M$ has a left inverse.\\
$(iii) \Rightarrow (ii)$. Let us denote by $\gamma: \gG\fF(M) \to M$
the left inverse of $\theta_M$. For any sequence
$$\xymatrix@C=30pt{ E: X \ar@{->}^-{i}[r] &
X' \ar@{->}^-{j}[r] & X''}$$ in $\eE$, we need to prove that its
corresponding sequence of abelian groups
\[
\xymatrix{ \Hom_\aA(X'',M) \ar[r] & \Hom_\aA(X',M) \ar[r] &
\Hom_\aA(X,M) }
\]
is exact (in the usual sense). Given such $E$ in $\eE$, we have a
commutative diagram in $\bB$ $$\xymatrix@C=30pt{ \fF(X)
\ar@{->}^-{\fF(i)}[r] & \fF(X') \ar@{->}^-{\fF(j)}[r]
\ar@{->}_-{\fF(i)^c}[d] & \fF(X'')
\\ & \mathrm{Coker}(\fF(i)) \ar@{->}_-{l}[ur]
& }$$ where $l$ splits as monomorphism by $l'$. Let $\tau: X' \to M$
be a morphism in $\aA$, such that $\tau \circ i = 0$. Then there
exists a morphism $g:\mathrm{Coker}(\fF(i)) \to \fF(M)$ of $\bB$
such that $g \circ \fF(i)^c \,=\, \fF(\tau)$. This leads to the
composition
$$\xymatrix@C=60pt{X'' \ar@{->}_-{\theta_{X''}}[d]
\ar@{-->}^-{\alpha}[r] & M \\ \gG \fF(X'') \ar@{->}^-{\gG(g \circ
l')}[r] & \gG \fF(M) \ar@{->}_-{\gamma}[u] }$$ The morphism
$\alpha$ satisfies
\begin{eqnarray*}
  \alpha \circ j &=& \gamma \circ \gG(g \circ l') \circ
\theta_{X''} \circ j \\
   &=& \gamma \circ \gG(g \circ l') \circ \gG\fF(j) \circ \theta_{X'} \\
   &=& \gamma \circ \gG\lr{ g \circ l' \circ \fF(j)} \circ \theta_{X'} \\
   &=& \gamma \circ \gG\lr{ g \circ l' \circ l \circ \fF(i)^c} \circ \theta_{X'} \\
   &=& \gamma \circ \gG\lr{ g \circ \fF(i)^c } \circ \theta_{X'} \\
   &=& \gamma \circ \gG\fF(\tau) \circ \theta_{X'} \\
   &=& \gamma \circ \theta_M \circ \tau \,\, =\,\, \tau
\end{eqnarray*} which proves the exactness of the sequence of abelian
groups.\\ $(ii) \Rightarrow (i)$ Let $i:X \to X'$ be a morphism of
$\aA$ such that $\fF(i)$ has a left inverse. The later condition
means that $\xymatrix{0 \ar@{->}[r]& \fF(X) \ar@{->}^-{\fF(i)}[r] &
\fF(X')}$ is a cosplit sequence in $\bB$. Thus $\xymatrix{0
\ar@{->}[r]& X \ar@{->}^-{i}[r] & X'}$ is a sequence in $\eE$.
Therefore, the corresponding sequence of abelian groups
$$\xymatrix{\hom{\aA}{X'}{M} \ar@{->}[r]& \hom{\aA}{X}{M}  \ar@{->}[r] &
0}$$ is exact. Whence $\hom{\aA}{i}{M}$ is surjective and so $M$ is
$\fF$-injective.
\end{proof}

Consider in the category of bicomodules $\lBicomod{\Aa}{\bf{F}}$
the class $\eE_0$ of all co-split sequences. This an injective
class, as $\lBicomod{\Aa}{\bf{F}}$ is an additive category with
cokernels. As we have mention, the corresponding class of
$\eE_0$-injective objects coincides with the class of all objects
of $\lBicomod{\Aa}{\bf{F}}$. Denote by $\eE := \sS^{-1}\lr{\eE_0}$
the class of sequences $E$ in the category
$\lrBicomod{\bf{F}}{\bf{F}}$ such that $\sS(E)$ is a sequence in
$\eE_0$, as we have point out $\eE$ is also an injective class.

\begin{prop}\label{S-ijectivos}
Let $(\Sf{M},\mmm,\nnn)$ be an $\bf{F}$-bicomodule. The following
are equivalent
\begin{enumerate}[(i)]
\item $(\Sf{M},\mmm,\nnn)$ is $\eE$-injective.

\item $(\Sf{M},\mmm,\nnn)$ is $\sS$-injective.

\item The unit $\Theta_{(\Sf{M},\,\mmm,\,\nnn)}$ of the adjunction
$\sS \dashv \tT$ at $(\Sf{M},\mmm,\nnn)$, stated in \eqref{Theta},
is a split monomorphism.
\end{enumerate}
In particular every bicomodule of the form $\tT(\Sf{N},\rrr)$ is
$\eE$-injective, for every bicomodule $(\Sf{N},\rrr) \in
\lBicomod{\Aa}{\bf F}$, and so is every induced $\bf F$-bicomodule
$\tT\sS(\Sf{M},\mmm,\nnn)\,=\,(\Sf{M}F,\mmm_F,\Sf{M}\delta)$.
\end{prop}
\begin{proof}
Follows immediate from Lemma \ref{lema-Sinj}.
\end{proof}

Fix ${\bf{F}}=(F,\delta,\xi)$ a comonad in a Grothendieck category
$\Aa$ with $F \in \ContFunt{\Aa}{\Aa}$. For every
$\bf{F}$-bicomodule $(\Sf{M},\mmm,\nnn)$ and each $i \geq 1$, we
consider the $i^{-th}$ induced $\bf{F}$-bicomodule
$(\Sf{M}F^i,\mmm_{F^i}, \Sf{M}F^{i-1}\delta)$.

\begin{prop}\label{resolution}
Let $(\Sf{M},\mmm,\nnn)$ be any $\bf{F}$-bicomodule. The following
sequence in the category of $\bf{F}$-bicomodules
\begin{equation}\label{sequence}
\xymatrix@C=20pt{0 \ar@{->}^-{}[r] & \Sf{M} \ar@{->}^-{\nnn}[r] &
\Sf{M}F \ar@{->}^-{\ddd^0}[r] & \Sf{M}F^2 \ar@{->}^-{\ddd^1}[r]
&\cdots \ar@{->}^-{}[r] & \Sf{M}F^{n+1} \ar@{->}^-{\ddd^n}[r] &
\Sf{M}F^{n+2} \ar@{->}^-{}[r]& \cdots }
\end{equation} where $\ddd^0\,=\, \Sf{M}\delta - \nnn_F$ and recursively
\begin{equation}\label{ddd}
\ddd^{n+1} \,\,\,=\,\,\, \ddd^n_F \,+\, (-1)^{n+1}
\Sf{M}F^{n+1}\delta,\quad n=0,1,2,\cdots
\end{equation}
defines an $\eE$-injective resolution for $(\Sf{M},\mmm,\nnn)$.
\end{prop}

\begin{proof}
Let us denote by $E(\Sf{M})$ the sequence defined in
\eqref{sequence}. One can easily check that the family of morphisms
$$\Sf{u}_n:=(-1)^{n+1}\Sf{M}F^n\xi:\Sf{M}F^{n+1}
\longrightarrow  \Sf{M}F^{n}$$ in $\lBicomod{\Aa}{\bf{F}}$, defines
a contracting homotopy for $\sS(E(\Sf{M}))$. This implies by
\cite[Lemma 2.4]{Jonah:1968} that $\sS(E(\Sf{M}))$ is sequence in
$\eE_0$. Hence $E(\Sf{M})$ is in $\eE$.
\end{proof}

Let $(\Sf{N},\rrr,\sss)$ be another $\bf{F}$-bicomodule and denote
by $\Ext{ }{N}{M}$ the homology of the complex
\begin{equation}\label{seq-1}
\xymatrix@C=20pt{ 0 \ar@{->}[r] &
\hom{\,\lrBicomod{\bf{F}}{\bf{F}}}{\Sf{N}}{\Sf{M}F}\ar@{->}[r] &
\hom{\,\lrBicomod{\bf{F}}{\bf{F}}}{\Sf{N}}{\Sf{M}F^2} \ar@{->}[r] &
\cdots }
\end{equation}
obtained by applying the functor
$\hom{\,\lrBicomod{\bf{F}}{\bf{F}}}{\Sf{N}}{-}$ to the
$\eE$-injective resolution of $\Sf{M}$ given in \eqref{sequence}.
Using the natural isomorphism stated in Lemma \ref{adj}, we show
that the complex \eqref{seq-1} is isomorphic to
\begin{equation}\label{seq-2}
\xymatrix@C=20pt{ 0 \ar@{->}[r] &
\hom{\,\lBicomod{\Aa}{\bf{F}}}{\Sf{N}}{\Sf{M}}\ar@{->}^-{\partial^0}[r]
& \hom{\,\lBicomod{\Aa}{\bf{F}}}{\Sf{N}}{\Sf{M}F}
\ar@{->}^-{\partial^1}[r] & \cdots }
\end{equation}
where
\begin{eqnarray*}
  \partial^0(\Sf{f}) &=& \Sf{f}_F \circ \sss - \nnn \circ \Sf{f}, \\
 \partial^1(\Sf{f})  &=& \Sf{M}\delta \circ \Sf{f} - \Sf{f}_F \circ \sss - \nnn_F\circ\Sf{f},   \\
 \partial^n(\Sf{f})  &=& \sum_{i=0}^{n-1} (-1)^i \Sf{M}F^i
 \delta_{F^{n-i-1}} \circ \Sf{f} + (-1)^n \Sf{f}_F \circ \sss -
 \nnn_{F^n} \circ \Sf{f},\quad n=2,3,...
\end{eqnarray*}
In particular, we have
\begin{eqnarray*}
  \Ker(\partial^1) &=& \left\{ \Sf{f}: (\Sf{N},\rrr) \to (\Sf{M}F,\mmm_F)|\,
  \Sf{M}\delta \circ \Sf{f} = \Sf{f}_F \circ \sss + \nnn_F\circ\Sf{f} \right\} \\
  \im(\partial^0) &=& \left\{ \Sf{f}:(\Sf{N},\rrr) \to (\Sf{M}F,\mmm_F)
  |\, \Sf{f}= \varphi_F \circ \sss - \nnn \circ \varphi,\,\, \text{ for}\,\,\text{ some }
  \,\, \varphi:(\Sf{N},\rrr) \to (\Sf{M},\mmm) \right\}
\end{eqnarray*}
That is the $1$-cocycle are cointegrations and the $1$-coboundaries
are inner cointegrations. Thus
\begin{equation}\label{Ext}
\Ext{1}{N}{M} \,\,\cong\,\, \Coint{\Sf{N}}{\Sf{M}} /
\IntCoint{\Sf{N}}{\Sf{M}}.
\end{equation}

The pair $(\tT\sS, \Theta_{-})$ form a resolvent pair in the sense
of \cite[Prposition 2.10]{Jonah:1968} for the injective class
$\eE$. Since $\lrBicomod{\bf{F}}{\bf{F}}$  has cokernels,
\cite[Lemma 2.11]{Jonah:1968} implies that the cokernels
constructed in \eqref{sucesion} lead to a functor $$\kK:
\lrBicomod{\bf{F}}{\bf{F}} \to \lrBicomod{\bf{F}}{\bf{F}},$$ and a
natural transformation
$$\tT\sS \to \kK.$$
Furthermore, $\kK(E)$ is a sequence in $\eE$, whenever $E$ is a
sequence in $\eE$. By the isomorphism given in \eqref{Isomofos},
we have $\hom{\,\lrBicomod{\bf{F}}{\bf{F}}}{\Sf{N}}{\kK(E)} \cong
\Coint{\Sf{N}}{E}$ is an exact sequence of abelian groups, for
every $\eE$-projective $\bf{F}$-bicomodule $\Sf{N}$ and every
sequence $E$ in $\eE$. On the other hand, given an $\eE$-injective
$\bf{F}$-bicomodule $\Sf{M}$, then $\kK(\Sf{M})$ is clearly
$\eE$-injective. Thus $\Coint{E}{\Sf{M}}$, which by
\eqref{Isomofos} is isomorphic to
$\hom{\,\lrBicomod{\bf{F}}{\bf{F}}}{E}{\kK(\Sf{M})}$, is an exact
sequence of abelian groups. This proves that the $\eE$-derived
functor of the bifunctor $\Coint{-}{-}$ can be constructed. For
$\Sf{N}$ and $\Sf{M}$ two $\bf{F}$-bicomodules, let
$H^*(\Sf{N},\Sf{M})$ be this $\eE$-derived functor which can be
computed using the $\eE$-injective resolution given in proposition
\ref{resolution}. Using this times the natural isomorphisms of
\eqref{Isomofos} and the fact that $\tT\sS(\Sf{M})$ are
$\eE$-injective for every $\bf{F}$-bicomodule $\Sf{M}$, we can
easily show that
\begin{eqnarray}
\Ext{n}{\Sf{N}}{\kK(\Sf{M})} &\cong & H^n(\Sf{N},\Sf{M}), \qquad n\geq 0  \label{Ext-1}\\
   \Ext{n+1}{\Sf{N}}{\Sf{M}} &\cong &
   \Ext{n}{\Sf{N}}{\kK(\Sf{M})},
   \qquad n \geq 1. \label{Ext-2}
\end{eqnarray}

By both propositions \ref{S-ijectivos} and \ref{univ-coint}, and
the isomorphisms given in  \eqref{Ext}, \eqref{Ext-1}, and
\eqref{Ext-2}, we have \begin{cor}\label{E-inj} For a
$\bf{F}$-bicomodule $(\Sf{M},\mmm,\nnn)$, the following are
equivalent
\begin{enumerate}[(i)]
\item $\Sf{M}$ is $\eE$-injective.

\item $\Sf{M}$ is $\sS$-injective.

\item The sequence $$\xymatrix@C=30pt{ 0 \ar@{->}[r] & \Sf{M}
\ar@{->}^-{\nnn}[r] & \Sf{M}F \ar@{->}^-{\nnn^c}[r] & \kK(\Sf{M})
\ar@{->}[r] & 0 }$$ splits in the category of bicomodules
$\lrBicomod{\bf{F}}{\bf{F}}$.

\item The universal cointegration from $\kK(\Sf{M})$ into $\Sf{M}$
is inner.

\item Every cointegration into $\Sf{M}$ is inner.
\end{enumerate}
\end{cor}

Now we can formulate a characterization of comonads with a
separable forgetful functor by means of the cohomology groups of
their bicomodules.

\begin{theo}\label{characterization}
Let $\Aa$ be a Grothendieck category and ${\bf{F}} = (F,\delta,\xi)$
a comonad in $\Aa$ with universal cogenerator the adjunction
$\xymatrix{ {\bf S}: \Aa^{\bf{F}} \ar@<0,5ex>[r] & \ar@<0,5ex>[l]
\Aa: \bf{T}}$. If $F$ is right exact and preserves direct sums. Then
the following are equivalent
\begin{enumerate}[(i)]
\item ${\bf S}:\Aa^{\bf{F}} \to \Aa$ is a separable functor.

\item $\sS:\lrBicomod{\bf{F}}{\bf{F}} \to \lBicomod{\Aa}{\bf{F}}$
is a Maschke functor.

\item $\delta: F \rightarrow F^2$ is a split monomorphism in the
category of bicomodules $\lrBicomod{\bf{F}}{\bf{F}}$.

\item $(F,\delta,\delta)$ is $\eE$-injective $\bf{F}$-bicomodule.

\item The universal coderivation from $\kK(F)$ into $F$ is inner.

\item Every coderivation into $F$ is inner.

\item All cointegrations between $\bf{F}$-bicomodules are inner.

\item $\Ext{n}{-}{-} \,=\,0$ for all $n \geq 1$.

\item $H^n(\Sf{N}, F)\,=\,0$ for all $\bf{F}$-bicomodule $\Sf{N}$
and all $n \geq 1$.
\end{enumerate}
\end{theo}
\begin{proof}
Corollary \ref{E-inj}, Proposition \ref{Coint-Coder}, and properties
of $\rm Ext$ give the following equivalences $ (ii) \Leftrightarrow
(vii)$, $(ii) \Leftrightarrow (viii)$, $(iv) \Leftrightarrow (ix)$,
$(iv) \Leftrightarrow (vi)$. Proposition \ref{S-ijectivos} gives the
equivalence $(iv) \Leftrightarrow (iii)$, and lastly Theorem
\ref{Separable} gives the equivalences $(i) \Leftrightarrow (ii)
\Leftrightarrow (iii)$.
\end{proof}

\section{Applications}

We present in this section two different applications of Theorem
\ref{characterization}. The first one is devoted to a coseparable
corings \cite{Guzman:1989}, where of course the comonad is defined
by the tensor product over algebra. The second deals with the
co-algebra coextension over fields, and the comonad is defined using
cotensor product. Here we obtain Nakajima's results \cite{N} without
requiring the co-commutativity of the base co-algebra. This
condition is however replaced, in our case, by assuming that the
extended coalgebra is a left co-flat.

\subsection{Coseparable corings}

Let $\mathbb{K}$ be commutative ring with $1$. In what follows all
algebras are $\mathbb{K}$-algebras, and all bimodules over algebras
are assumed to be central $\mathbb{K}$-bimodules. Let $R$ be an
algebra an \emph{$R$-coring} \cite{Sweedler:1975} is a three-tuple
$(\coring{C}, \Delta, \varepsilon)$ consisting of a $R$-bimodule and
two $R$-bilinear maps $$\Delta: \coring{C} \to \coring{C} \tensor{R}
\coring{C}\quad \textrm{and}\quad \varepsilon: \coring{C} \to R,$$
known as the comultiplication and the counit, which satisfy
$$(\coring{C}\tensor{R}\Delta) \circ \Delta \,\, =\,\, (\Delta
\tensor{R} \coring{C}) \circ \Delta,\quad (\coring{C} \tensor{R}
\varepsilon) \circ \Delta \,\,=\,\, \coring{C} \,\,=\,\,
(\varepsilon\tensor{R}\coring{C}) \circ \Delta.$$ In this
sub-section the unadorned symbol $-\tensor{}-$ between
$R$-bimodules and $R$-bilinear maps denotes the tensor product
$-\tensor{R}-$. We denote as usual by
$\Bicomod{\coring{C}}{\coring{C}}$ the category of
$\coring{C}$-bicomodules. The objects are three-tuples
$(M,\varrho_M,\lambda_M)$ consisting of $R$-bimodule $M$ and two
$R$-bilinear maps $\varrho_M: M \to M \otimes\coring{C}$ (right
$\coring{C}$-coaction), $\lambda_M: M \to \coring{C}\otimes M$
(left $\coring{C}$-coaction) satisfying
\begin{eqnarray*}
  (\coring{C}\tensor{}
\lambda_M) \circ \lambda_M &=& (\Delta \tensor{} M) \circ
\lambda_M,\quad (\varepsilon \otimes M) \circ \lambda_M \,\,=\,\, M \\
  (\varrho_M \tensor{} \coring{C}) \circ \varrho_M &=& (M\otimes \Delta) \circ \varrho_M,
  \quad (M\otimes \varepsilon) \circ \varrho_M \,\,=\,\,M  \\
 (\coring{C}\otimes \varrho_M) \circ \lambda_M  &=& (\lambda_M
 \otimes \coring{C}) \circ \varrho_M.
\end{eqnarray*}

It is clear that ${\bf F}:=(F,\delta,\xi)$ where $F=-\otimes
\coring{C}: \rmod{R} \to \rmod{R}$, $\delta=-\otimes \Delta$, and
$\xi=-\otimes \varepsilon$, is a comonad in the category of right
$R$-modules $\rmod{R}$, with $F \in
\ContFunt{\rmod{R}}{\rmod{R}}$.

Given any $\bf{F}$-bicomodule $(\Sf{M},\mmm,\nnn)$ we can use
Watts' theorem \cite{Watts:1960} to find a natural isomorphism
\begin{equation}\label{daleth}
\daleth_{-}^{\Sf{M}}: \Sf{M} \longrightarrow -\otimes \Sf{M}(R)
\end{equation}
satisfying $(-\otimes \psi_R) \circ \daleth_{-}^{\Sf{M}} \,\, =\,\,
\daleth_{-}^{\Sf{M}'} \circ \psi$ for every natural transformation
$\psi:\Sf{M} \to \Sf{M}'$ with $(\Sf{M}',\mmm',\nnn')$ is another
$\bf{F}$-bicomodule. With the help of this natural isomorphism we
can establish a functor
$$\xymatrix@C=50pt@R=0pt{ \gG:\,\, \lrBicomod{\bf{F}}{\bf{F}} \ar@{->}[r] &
\Bicomod{\coring{C}}{\coring{C}} \\ (\Sf{M},\mmm,\nnn) \ar@{->}[r]
& \lr{\Sf{M}(R),\varrho_{\Sf{M}(R)},\lambda_{\Sf{M}(R)}} \\
\Sf{f} \ar@{->}[r] & \Sf{f}_R}$$ where the $\coring{C}$-coactions
are defined by $\varrho_{\Sf{M}(R)}\,=\, \mmm_R$ and
$\lambda_{\Sf{M}(R)}\,=\, \daleth_{F(R)}^{\Sf{M}} \circ \nnn_R$.

Conversely, given any $\coring{C}$-bicomodule
$(M,\varrho_M,\lambda_M)$, we clearly obtain a $\bf{F}$-bicomodule
defined by the three-tuple $\lr{-\otimes M, -\otimes \varrho_M,
\lr{\daleth_{F}^{\Sf{M}}}^{-1} \circ (-\otimes \lambda_M)}$. This
in fact entails an inverse functor, up to the natural isomorphisms
$\daleth_{-}^{-}$, to the functor $\gG$. Henceforth, $\gG$ is an
equivalence of categories $\lrBicomod{\bf{F}}{\bf{F}}$ and
$\Bicomod{\coring{C}}{\coring{C}}$. It is then obvious that
$\delta$ is a split-mono in the category of $\bf{F}$-bicomodules
if and only if $\Delta$ is a split-mono in the category of
$\coring{C}$-bicomodules. It is well known (see
\cite{Brzezinski/Wisbauer:2003}) that this later condition happens
if and only if the right coaction forgetful functor is separable.

Given two $\coring{C}$-bicomodules $(M,\varrho_M,\lambda_M)$ and
$(N,\varrho_N,\lambda_N)$. Following to \cite{Guzman:1989}, a
$R$-bilinear map $g: M \to \coring{C}$ is said to be
\emph{coderivation} if it satisfies
$$ \Delta \circ g \,\,=\,\, (g \otimes \coring{C}) \circ \varrho_M
+ (\coring{C}\otimes g) \circ \lambda_M$$ The coderivation $g$ is
said to be an \emph{inner coderivation} if there exists a
$R$-bilinear map $\gamma:M \to R$ such that $g =
(\coring{C}\otimes \gamma) \circ \lambda_M - (\gamma \otimes
\coring{C}) \circ \varrho_M$. We denote by ${\rm
Coder}_{\coring{C}}(M,A)$ the abelian group of all coderivations
from $M$ to $\coring{C}$. A \emph{(left) cointegration} from $N$
into $M$ is a $R$-bilinear morphism $f: N \to \coring{C} \otimes
M$ such that
$$(\Delta \otimes \coring{C}) \circ f \,\, =\,\, (\coring{C} \otimes
\lambda_M) \circ f + (\coring{C}\otimes f) \circ \lambda_N$$ The
cointegration $f$ is said to be an \emph{inner cointegration} if
there exists a $R$-bilinear map $\varphi: N \to M$ satisfying
$$ \varrho_M \circ \varphi \,\, =\,\, (\varphi \otimes \coring{C}) \circ
\varrho_N,\,\, \text{and} \,\, f\,\,=\,\, (\coring{C}\otimes
\varphi) \circ \lambda_N - \lambda_M \circ \varphi$$ The abelian
group of all cointegrations from $N$ into $M$ will be denoted by
${\rm Coint}_{\coring{C}}(N,M)$.

Cointegrations and coderivations in both categories of bicomodules
$\lrBicomod{\bf{F}}{\bf{F}}$ and $\Bicomod{\coring{C}}{\coring{C}}$
are connected by the following isomorphisms of an abelian groups
$$\xymatrix@R=0pt@C=40pt{\Coder{\Sf{M}}{F} \ar@{->}^-{\cong}[r] &
{\rm Coder}_{\coring{C}}\lr{\Sf{M}(R),\coring{C}} \\
\Sf{g} \ar@{|->}[r] & \Sf{g}_R \\ (-\otimes g) \circ
\daleth_{-}^{\Sf{M}} & \ar@{|->}[l] g }$$ and
$$\xymatrix@R=0pt@C=40pt{\Coint{\Sf{N}}{\Sf{M}} \ar@{->}^-{\cong}[r] &
{\rm Coint}_{\coring{C}}\lr{\Sf{N}(R),\Sf{M}(R)} \\
\Sf{f} \ar@{|->}[r] &  \daleth_{F(R)}^{\Sf{M}}  \circ \Sf{f}_R \\
\lr{\daleth_{F(-)}^{\Sf{M}}}^{-1} \circ (-\otimes f) \circ
\daleth_{-}^{\Sf{N}} & \ar@{|->}[l] f }$$ where the isomorphism
$F(R) \cong \coring{C}$ was used as isomorphism of $R$-corings.
The restrictions of the above isomorphisms to the sub-groups of
inner coderivations or inner cointegrations, are also
isomorphisms.

Applying Theorem \ref{characterization} to this situation, we obtain
\begin{cor}[{\cite[Theorem 3.10]{Guzman:1989}}]
For any $R$-coring $(\coring{C},\Delta,\varepsilon)$, the
following are equivalent
\begin{enumerate}[(i)]
\item The forgetful functor ${\bf S}:\rcomod{\coring{C}} \to
\rmod{R}$ from the category of right $\coring{C}$-comodules to the
category of right $R$-modules is a separable functor.

\item The forgetful functor ${^\cc\mM^\cc}\to {^\cc\mM_R}$ is a
Maschke functor.

\item The short exact sequence $$\xymatrix@C=40pt{ 0 \ar@{->}[r] &
\coring{C} \ar@{->}^-{\Delta}[r] & \coring{C}\otimes \coring{C}
\ar@{->}^-{\Delta^c}[r] & \Omega(\coring{C}) \ar@{->}[r] & 0 }$$
splits in the category of bicomodules
$\Bicomod{\coring{C}}{\coring{C}}$.

\item $\coring{C}$ is $\eE$-injective, where $\eE$ is the
injective class in $\Bicomod{\coring{C}}{\coring{C}}$ whose
sequences split in the category of $R$-bimodules $\Bimod{R}{R}$.

\item The universal coderivation from $\Omega(\coring{C})$ into
$\coring{C}$ is inner.

\item Every coderivation into $\coring{C}$ is inner.

\item All cointegrations between $\coring{C}$-bicomodules are
inner.

\item $\Ext{n}{-}{-} \,=\,0$ for all $n \geq 1$.

\item $H^n(N, \coring{C})\,=\,0$ for all $\coring{C}$-bicomodule
$N$ and all $n \geq 1$.
\end{enumerate}
\end{cor}

\subsection{Coseparable coalgebras co-extension}
In what follows $\mathbb{K}$ is assumed to be a field. The unadorned
symbol $\otimes$ between $\mathbb{K}$-vector spaces means the tensor
product $\otimes_{\mathbb{K}}$. Let $A$, $C$ are two
$\mathbb{K}$-coalgebras, and consider $\phi: A \to C$ a morphism of
$\mathbb{K}$-coalgebras. This define an adjunction
$$\xymatrix{ -\square_C A: \rcomod{C} \ar@<0,5ex>[r] & \ar@<0,5ex>[l] \rcomod{A}:
\oO}$$ between the categories of right comodules with
$-\square_CA$ right adjoint to $\oO$, and where $-\square_C-$ is
the co-tensor product over $C$. In the remainder, we denote this
bi-functor by $-\square-:=-\square_C-$. Notice that $-\square-$ is
associative (up to natural isomorphism), as $C$ is a
$\mathbb{K}$-coalgebra and $\mathbb{K}$ is a field. From now on,
we assume that $-\square A: \rcomod{C} \to \rcomod{A}$ is right
exact, and thus exact. Put $F:=\oO (-\square A): \rcomod{C} \to
\rcomod{C}$, since $\rcomod{C}$ is a Grothendieck category we can
construct the category $\ContFunt{\rcomod{C}}{\rcomod{C}}$, and we
have in this case that $F \in \ContFunt{\rcomod{C}}{\rcomod{C}}$.
Let us denote by $\bara{\Delta}: A \to A\square A$ the resulting
map from the universal property of kernels. This is in fact an
$A$-bicomodule map, and thus a $C$-bicomodule map by applying
$\oO$. Furthermore, we have
\begin{eqnarray*}
(A\square \bara{\Delta}) \circ \bara{\Delta} &=& (\bara{\Delta} \square A) \circ \bara{\Delta} \\
  (\phi \square A) \circ \bara{\Delta} &=& (A \square \phi) \circ
  \bara{\Delta} \,\, =\,\, A \,\, \text{(up} \, \text{to}\,
  \text{isomorphisms)}.
\end{eqnarray*}
Using these equalities, on can easily check that there is a comonad
${\bf F}:=(F,\delta, \xi)$ in the category of right $C$-comodules
$\rcomod{C}$, where $\delta$ and $\xi$ are defined by the following
commutative diagrams of natural transformations
$$\xymatrix@C=40pt{ -\square A \ar@{->}^-{-\square\bara{\Delta}}[r] & -\square A\square A,
\\ F \ar@{=}[u] \ar@{-->}^-{\delta}[r] &  F^2 \ar@{=}[u]}
\qquad \xymatrix@C=40pt{-\square A \ar@{->}^-{-\square\phi}[r] &
-\square C \ar@{->}^-{\cong}[d] \\ F \ar@{=}[u]
\ar@{-->}^-{\xi}[r] & \id_{\rcomod{C}}}$$

Given $(\Sf{M},\mmm,\nnn)$ any ${\bf F}$-bicomodule, we know that
$\Sf{M}: \rcomod{C} \to \rcomod{C}$ is right exact and preserves
direct sums. If $\Sf{M}$ is assumed to be left exact , then by
\cite[Theorem 2.6]{Gomez-Torrecillas:2002}, $\Sf{M}(C):=M$ is a
$C$-bicomodule, and  there is a natural isomorphism
\begin{equation}\label{upsilon}
\xymatrix@C=50pt{\Upsilon^{\Sf{M}}_{-}\,:\, \Sf{M}
\ar@{->}^-{\cong}[r] & -\square M, }
\end{equation} which satisfies $(-\square\beta_C) \circ
\Upsilon_{-}^{\Sf{M}} \,\, =\,\, \Upsilon_{-}^{\Sf{N}} \circ
\beta$, for every natural transformation $\beta: \Sf{M} \to
\Sf{N}$ with $\Sf{N} \in \ContFunt{\rcomod{C}}{\rcomod{C}}$ and
$\Sf{N}$ an exact functor. The natural transformation $\mmm$ and
$\nnn$ induces by this isomorphism a structure of $A$-bicomodule
on $M$. The right and left $A$-coactions are given by
$$\xymatrix@C=50pt{M
\ar@{->}^-{\mmm_C}[r] \ar@{-->}_-{\varrho_M}[dr] & M\square A \ar@{->}^-{\equalizerk{M}{A}}[d] \\
& M \otimes A  }\qquad  \xymatrix@C=50pt{M \ar@{->}^-{\nnn_C}[r]
\ar@{-->}_-{\lambda_M}[drr] & \Sf{M}F(C)
\ar@{->}^-{\Upsilon^{\Sf{M}}_{F(C)}}_-{\cong}[r] & F(C)\square M
\cong A\square M
\ar@{->}^-{\equalizerk{A}{M}}[d]  \\
&  & A \otimes M  }$$ where $\equalizerk{X}{Y}$ is the equalizer
map, that is the kernel of $\xymatrix{\equalizer{X}{Y}: X \otimes
Y \ar@<0,7ex>[r] & \ar@{<-}[l]  X \otimes C \otimes Y}$ defined by
$\equalizer{X}{Y}\,=\, \varrho_X\tensor{}Y - X\tensor{}\lambda_Y$
for every right $C$-comodule $X$ and left $C$-comodule $Y$. The
counitary conditions of these new $A$-coactions are easily seen,
while the co-associatively and compatibility conditions need a
routine and long computations using properties of cotensor product
over coalgebras over fields. Let us denote by ${}^{\bf
F}\Sf{E}^{\bf F}$ the full subcategory of $\lrBicomod{\bf F}{\bf
F}$ whose objects are $\bf F$-bicomodules $(\Sf{M},\mmm,\nnn)$
such that $\Sf{M}: \rcomod{C} \to \rcomod{C}$ is an exact functor
which commutes with direct sums. For instance $(F,\delta,\xi)$ is
an object of this category.

The above arguments establishes in fact a functor from the
subcategory of $\bf{F}$-bicomodules ${}^{\bf F}\Sf{E}^{\bf F}$ to
the category of $A$-bicomodule sending
\begin{equation}\label{fF}
\fF: {}^{\bf F}\Sf{E}^{\bf F} \longrightarrow \Bicomod{A}{A},\,\,
\lr{\lr{\Sf{M},\mmm,\nnn} \to \lr{M,\varrho_{M},\lambda_{M}}},\,
\lr{ \Sf{f} \to \Sf{f}_C}
\end{equation}
For every $\bf{F}$-bicomodule $\Sf{M} \in {}^{\bf F}\Sf{E}^{\bf
F}$, it is clear that $\fF(\Sf{M})=M$ is a co-flat left
$C$-comodule.

Conversely, given any $A$-bicomodule $(N,\varrho_N,\lambda_N)$
such that the underlying left $C$-comodule ${}_CN$ is co-flat,
then we have a functor $-\square N: \rcomod{C} \to \rcomod{C}$
which is exact and preserves direct sums together with two natural
transformations
$$\xymatrix@C=50pt{-\square N \ar@{->}^-{-\square \lambda'_N}[r] & -\square A \square N, }
\qquad \xymatrix@C=50pt{-\square N \ar@{->}^-{-\square
\varrho'_N}[r] & -\square N \square A, }$$ where $\lambda_N'$ and
$\varrho_N'$ are $C$-bicolinear defined by universal property
$$\xymatrix@C=50pt{ N
\ar@{->}^-{\lambda_N}[r] \ar@{-->}_-{\lambda_N'}[rd] & A \otimes N \\
& A \square N \ar@{->}_-{\equalizerk{A}{N}}[u]} \qquad
\xymatrix@C=50pt{ N \ar@{->}^-{\varrho_N}[r] \ar@{-->}_-{\varrho_N'}[rd] & N \otimes A \\
& N \square A \ar@{->}_-{\equalizerk{N}{A}}[u]}$$ By definition and
the properties of cotensor product $\lambda'_N$ and $\varrho'_N$
satisfy the following equalities
\begin{eqnarray*}
  (\bara{\Delta} \square N) \circ \lambda_N' &=& (A \square \lambda_N') \circ \lambda_N',
  \quad (\phi \square N) \circ \lambda_N' \,\, =\,\, N \,\, \text{(up} \, \text{to}\,
  \text{isomorphism)}\\
  (N\square \bara{\Delta}) \circ \varrho_N '&=& (\varrho_N'
  \square A) \circ \varrho_N ',\quad (N\square \phi) \circ
  \varrho_N' \,\, =\,\, N \,\, \text{(up} \, \text{to}\,
  \text{isomorphism)}\\ (A\square \varrho_N') \circ
  \lambda_N' &=& (\lambda_N' \square A) \circ \varrho_N'.
\end{eqnarray*}
Consider the obtained three-tuple $(\Sf{N},\rrr,\sss)$, where
$\Sf{N}:=-\square N: \rcomod{C} \to \rcomod{C}$ is a functor, and
$\rrr:= -\square \varrho_N': \Sf{N} \to F\Sf{N}$, $\sss:=-\square
\lambda_N': \Sf{N} \to \Sf{N}F$ are two natural transformation.
Since $N$ is assumed to be co-flat left $C$-comodule, the previous
equalities show that $(\Sf{N},\rrr,\sss)$ is actually an object of
the category ${}^{\bf F}\Sf{E}^{\bf F}$, whose image by $\fF$ is
isomorphic to the initial $A$-bicomodule
$(N,\varrho_N,\lambda_N)$, via the natural isomorphisms
$\Upsilon_{-}^{-}$. Now, given an $A$-bicolinear morphism
$g:(N,\varrho_N,\lambda_N) \to (N',\varrho_{N'},\lambda_{N'})$, we
get a $\bf{F}$-bicomodules morphism $\Sf{g}:=-\square g: \Sf{N}
\to \Sf{N}'$. This shows that the above constructions are in fact
functorial.

In conclusion, we have shown that the functor $\fF$ defined in
\eqref{fF}, establishes an equivalence of categories ${}^{\bf
F}\Sf{E}^{\bf F}$ and ${}^A\cC^A$, where the later is the full
subcategory of the category of $A$-bicomodules $\Bicomod{A}{A}$
whose objects are co-flat left $C$-comodules after forgetting the
right $C$-coaction.

Recall from \cite{N} that $A$ is said  to be a \emph{separable
$C$-coalgebra} if the $A$-bicolinear map $\bara{\Delta}: A \to
A\square A$ is a split-mono in the category of $A$-bicomodules. By
\cite[Theorem 4.7]{Gomez-Torrecillas:2002} this is equivalent to
say that the forgetful functor $\oO$ is a separable functor. Using
the equivalence of categories established above, it is easy to
check that $\delta$ is a split-mono in ${}^{\bf F}\Sf{E}^{\bf F}$
(or equivalently in $\lrBicomod{\bf{F}}{\bf{F}}$) if and only if
$\bara{\Delta}$ is a split-mono in ${}^{A}\cC^A$ (or equivalently
in $\Bicomod{A}{A}$).

Given two $A$-bicomodules $(M,\varrho_M,\lambda_M)$ and
$(N,\varrho_N,\lambda_N)$, a $C$-bicolinear map $g: M \to A$ is
said to be \emph{$C$-coderivation} if its satisfies $$
\bara{\Delta} \circ g \,\,=\,\, (g \square A) \circ \varrho_M' +
(A\square g) \circ \lambda_M'$$ The $C$-coderivation $g$ is said
to be an \emph{inner $C$-coderivation} if there exists a
$C$-bicolinear map $\gamma:M \to C$ such that $g = (A\square
\gamma) \circ \lambda_M' - (\gamma \square A) \circ \varrho_M'$.
We denote by ${\rm Coder}_C(M,A)$ the abelian group of all
$C$-derivations from $M$ to $A$. A \emph{(left) $C$-cointegration}
from $N$ into $M$ is a morphism of $C-A$-bicomodule $f: N \to A
\square M$ such that
$$(\bara{\Delta} \square A) \circ f \,\, =\,\, (A \square
\lambda_M') \circ f + (A\square f) \circ \lambda_N'$$ The
$C$-cointegration $f$ is said to be an \emph{inner
$C$-cointegration} if there exists a $C$-bicolinear map $\varphi:
N \to M$ satisfying $$ \varrho_M' \circ \varphi \,\, =\,\,
(\varphi \square A) \circ \varrho_N',\,\, \text{and} \,\,
f\,\,=\,\, (A\square \varphi) \circ \lambda_N' - \lambda_M' \circ
\varphi$$ The abelian group of all $C$-cointegration from $N$ into
$M$ will be denoted by ${\rm Coint}_C(N,M)$.

Given $(\Sf{M},\mmm,\nnn)$ and $(\Sf{N},\rrr,\sss)$ two
$\bf{F}$-bicomodules in ${}^{\bf F}\Sf{E}^{\bf F}$ and consider
their associated $A$-bicomodule via the above equivalence of
categories $\fF$:
$$(\Sf{M}(C):=M,\varrho_M,\lambda_M) \quad \textrm{and}\quad
(\Sf{N}(C):=N,\varrho_N',\lambda_N').$$ We have an abelian group
isomorphism $$\xymatrix@R=0pt@C=40pt{\Coder{\Sf{M}}{F} \ar@{->}^-{\cong}[r] &  {\rm Coder}_C\lr{\Sf{M}(C),A} \\
\Sf{g} \ar@{|->}[r] & \iota_A \circ \Sf{g}_C \\ (-\square g) \circ
\Upsilon_{-}^{\Sf{M}} & \ar@{|->}[l] g }$$ where $\iota_{-}:
C\square- \to \id_{\rcomod{C}}$ is the obvious natural isomorphism.
The isomorphism of cointegrations groups is given by
$$\xymatrix@R=0pt@C=40pt{\Coint{\Sf{N}}{\Sf{M}} \ar@{->}^-{\cong}[r] &  {\rm Coint}_C\lr{\Sf{N}(C),\Sf{M}(C)} \\
\Sf{f} \ar@{|->}[r] & (\iota_A \square \Sf{M}(C)) \circ \Upsilon_{F(C)}^{\Sf{M}}  \circ \Sf{f}_C \\
\lr{\Upsilon_{F(-)}^{\Sf{M}}}^{-1} \circ (-\square f) \circ
\Upsilon_{-}^{\Sf{N}} & \ar@{|->}[l] f }$$

Of course the restrictions of those isomorphisms to the sub-groups
of inner cointegrations or inner coderivations are also groups
isomorphisms. Applying now theorem \ref{characterization}, we
arrive to the following

\begin{cor}[{compare with \cite[Theorem 1.2]{N}}]
Let $\phi :A \to C$ be a morphism of $\mathbb{K}$-coalgebras over a
field $\mathbb{K}$. Assume that ${}_CA$ is a co-flat left
$C$-comodule. The following are equivalent
\begin{enumerate}[(i)]

\item $A$ is a separable $C$-coalgebra.

\item For any $A$-bicomodule $M$ such that ${}_CM$ is co-flat,
every $C$-coderivation from $M$ to $A$ is inner.

\item For any pair of $A$-bicomodules $M$ and $N$ such that
${}_CM$ and ${}_CN$ are co-flat, every $C$-cointegration from $M$
into $N$ is inner.

\end{enumerate}

\end{cor}

\subsubsection*{Acknowledgement}
The first author would like to express his thanks to S. Caenepeel
for inviting him to the Vrije Universiteit Brussel in August 2004,
and all members of Mathematics Department for a very warm
hospitality. This author is supported by grant MTM2004-01406 from
the Ministerio de Educaci\'{o}n y Ciencia of Spain.

\providecommand{\bysame}{\leavevmode\hbox
to3em{\hrulefill}\thinspace}
\providecommand{\MR}{\relax\ifhmode\unskip\space\fi MR }
\providecommand{\MRhref}[2]{
} \providecommand{\href}[2]{#2}

\end{document}